\documentclass[12pt]{article}
\usepackage[english]{babel}
\usepackage{amssymb, amsmath, amsthm, verbatim}
\newcounter{parag}
\newcommand{\sect}[1]
{\refstepcounter{parag}
\begin{center} { \bf\S\,\theparag. #1} \end{center}}
\newtheorem{theorem}{Theorem}
\newtheorem{cor}{Corollary}
\newtheorem{lemma}{Lemma}[parag]
\newtheorem{prop}[lemma]{Proposition}

\theoremstyle{definition}
\newtheorem*{prf}{Proof}

\begin{document}

\noindent{UDC 512.542}

\begin{center}\textbf{Spectra of finite symplectic and orthogonal groups}\end{center}

\begin{center}\textsc{Buturlakin A.A. \footnote{The work is supported by the Russian Foundation for
Basic Research (project 08--01--00322), the Council for Grants (under RF President) and State Aid
of Leading Scientific Schools (project NSh--3669.2010.1), the Program ``Development of the
Scientific Potential of Higher School'' of the Russian Federal Agency for Education
(project~2.1.1.419), the Federal Program ``Scientific and Scientific-Pedagogical Personnel of
Innovative Russia'' in 2009-2013 (gov. contract no. 02.740.11.0429), and also by Lavrent'ev Young
Scientists Competition (No 43 on 04.02.2010).}}
\end{center}

\begin{center}\textbf{Introduction}\end{center}

\textit{The spectrum} $\omega(G)$ of a finite group $G$ is the set of its element orders. In the
present paper we describe the spectra of groups $Sp_{2n}(q)$, $PSp_{2n}(q)$, $\Omega_{2n+1}(q)$,
$\Omega^\pm_{2n}(q)$, $P\Omega^\pm_{2n}(q)$, and also of groups $SO_{2n+1}(q)$ and $SO^\pm_{2n}(q)$
for odd $q$. In particular, a description of spectra of all finite simple simplectic and orthogonal
groups is obtained.

Let $G$ be a finite group of Lie type over a field of characteristic $p$. The set $\omega(G)$ can be presented as a
union of three subsets: the subset $\omega_p(G)$ of orders of all \textit{unipotent} elements, i.e. those whose order
is a power of $p$, the subset $\omega_{p'}(G)$ of orders of all \textit{semisimple} elements, i.e. those whose order is
coprime with $p$, and the subset $\omega_m(G)$ of the rest "composite"\ orders. Thus, the problem of describing the
spectrum of a finite group of Lie type splits into three subproblems.

The maximal order of unipotent elements in every finite group of Lie type is found in \cite[Proposition 0.5]{Tes}. It
is well known that each semisimple element of a group of Lie type is contained in some maximal torus of this group. In
\cite{our} a description of cyclic structure of maximal tori in all groups under consideration was obtained, and thus
the semisimple parts of spectra of these groups was described. Hence, it remains to describe the composite parts of the
spectra.

The set $\omega(G)$ is closed under taking divisors, i.e. if $n\in\omega(G)$ and $d$ divides $n$
then $d\in\omega(G)$. Therefore, it is uniquely determined by its subset $\mu(G)$ of elements that
are maximal under the divisibility relation. Define a set $\mu_m(G)$ to be the intersection of
$\mu(G)$ and $\omega_m(G)$. In the present paper for each group $G$ under consideration we
construct a set $\nu(G)$ such that $\mu_m(G)\subseteq\nu(G)\subseteq\omega(G)$.

The author would like to thank Vasil'ev A.V., Grechkoseeva M.A. and Vdovin E.P. for helpful remarks
on the paper.

\sect{Preliminary information and results}

The information on the theory of algebraic groups which we use in the following discussion can be
found, for example, in \cite[Chapters 1 and 3]{CarterBook}.

Let $p$ be a prime and $\overline{G}$ be a connected reductive algebraic group over algebraic
closure $\overline F_p$ of the field $GF(p)$. Surjective endomorphism $\sigma$ of $\overline G$ is
called a Frobenius map, if the fixed point group $\overline G_\sigma=\{g\in\overline G|
g^\sigma=g\}$ is finite. In this situation the finite group $\overline G_\sigma$ is called a finite
group of Lie type. А maximal closed connected diagonalizable subgroup of an algebraic group is
called a maximal torus. A subgroup of an algebraic group is called reductive, if its unipotent
radical is equal to identity. A reductive subgroup is called a reductive subgroup of maximal rank,
if it contains some maximal torus. \textit{A reductive subgroup of maximal rank} of
$\overline{G}_\sigma$ is a subgroup of the form $(\overline G_1)_\sigma$, where $\overline G_1$ is
a $\sigma$-stable reductive subgroup of maximal rank of $\overline{G}$.

If $g\in \overline G_\sigma$ then $g$ can be uniquely presented as a product $g_pg_{p'}$, where $g_{p'}$ is a
semisimple element and $g_{p}$ is a unipotent element of $C_{\overline G_\sigma}(g_{p'})$. Thus to describe the
composite part of spectrum it is sufficient for each semisimple element of $\overline G_\sigma$ to find the maximal
element of $\omega_p(C_{\overline G_\sigma}(s))$. In solving this problem, we use an approach developed in
\cite{Carter2, Carter}. Now we will formulate results from these works that we will need.

Let $s\in \overline{G}_\sigma$ be a semisimple element. Denote by $C_{\overline{G}}(s)^0$ the
connected component of centralizer $C_{\overline{G}}(s)$ containing identity. Subgroup
$C_{\overline{G}}(s)^0$ is a reductive subgroup of maximal rank of $\overline G$. Besides, subgroup
$C_{\overline{G}}(s)^0$ contains $s$ and all unipotent elements of $C_{\overline{G}}(s)$.
Therefore, reductive subgroup $(C_{\overline{G}}(s)^0)_\sigma$ of $\overline{G}_\sigma$ contains
$s$ and all unipotent elements of $C_{\overline{G}_\sigma}(s)$. Thus, the problem is reduced to the
following one: for every reductive subgroup of maximal rank of $\overline{G}_\sigma$ to find the
period of its center and the unipotent part of spectrum. For a reductive subgroup $G_1$, we will
denote by $\eta(G_1)$ the product of the period of its center and the maximal element of
$\omega_p(G_1)$.

Let $\overline T$ be a maximal torus of $\overline G$ and $\Phi$ be a root system of $\overline{G}$
with respect to $\overline T$. Subset $\Phi_1$ of $\Phi$ is called \emph{a subsystem} if it is a
root system itself. Subsystem is called \emph{closed} if for every $r_1, r_2 \in \Phi_1$ inclusion
$r_1+r_2\in \Phi$ yields inclusion $r_1+r_2\in \Phi_1$. For root $r\in \Phi$, denote by  $\overline
U_r$ the root subgroup corresponding to $r$. If $\overline G$ is a classical group distinct from
symplectic group over a filed of characteristic $2$, then reductive subgroups containing $\overline
T$ have form $\langle \overline T,\overline U_r, r\in\Phi_1\rangle$, where $\Phi_1$ is a closed
subsystem of $\Phi$. Root systems of reductive subgroups of $Sp_{2n}(\overline F_2)$ can be not
closed subsystems of its root system. However, group $Sp_{2n}(\overline F_2)$ is isomorphic to
$\Omega_{2n+1}(\overline F_2)$ and these symplectic groups can be considered as orthogonal groups.
Hence, in all cases witch are considered in the paper we can assume that root systems of reductive
subgroups are closed subsystems of $\Phi$. Henceforth, when we say "subsystem"\ we mean closed
subsystem.

Let $\overline G_1$ be a $\sigma$-stable reductive subgroup of  $\overline G$. Then $\overline G_1$ contains a
$\sigma$-stable maximal torus $\overline T$. Let $\overline G_1^g$ be also $\sigma$-stable. Then $\overline G_1^g$
contains a $\sigma$\nobreakdash-stable maximal torus. Since any two maximal tori of $\overline G_1^g$ are conjugate,
one can assume that torus $\overline T^g$ is $\sigma$-stable. Then $g^\sigma g^{-1}\in N_{\overline{G}}(\overline
T)\cap N_{\overline{G}}(\overline G_1)$. Let $W$ be the Weyl group of $\overline{G}$ and $W_1$ be the Weyl group of
$\overline{G}_1$. Denote by $\pi$ the canonical homomorphism from $N_{\overline{G}}(\overline T)$ to $W$. Then? as it
is shown in the proof of Proposition 2 from \cite{Carter2}, $\pi(N_{\overline{G}}(\overline T)\cap
N_{\overline{G}}(\overline G_1))=N_W(W_1)$.

Since $\overline{T}$ is stable under the action of $\sigma$, map $\sigma$ acts on $W=N_{\overline{G}}(\overline
T)/\overline T$. Since $W_1$ is also $\sigma$-stable, $\sigma$ acts on $N_W(W_1)/W_1$. Two elements $W_1w_1$ and
$W_1w_2$ of $N_W(W_1)/W_1$ are called \textit{$\sigma$-conjugate}, if there is $w\in W$ such that
$W_1w_2=(W_1w)^\sigma(W_1w_1)(W_1w)^{-1}$.

\begin{lemma}{\rm \cite[Proposition 3]{Carter2}}\label{l:bijection} Let $\overline G$ be a connected reductive
algebraic group and $\sigma$ be a Frobenius map of $\overline G$. Let $\overline G_1$ be a $\sigma$-stable reductive
subgroup of maximal rank of $\overline G$. Let $\mathfrak{C}$ be the set of all $\sigma$-stable subgroups conjugate to
$\overline G_1$ in $\overline G$ and  $\mathfrak{C}/\overline G_\sigma$ be the set of $\overline G_\sigma$-orbits in
$\mathfrak{C}$. The map $\pi$ induces a bijection between $\mathfrak{C}/\overline G_\sigma$ and classes of
$\sigma$-conjugate elements of $N_W(W_1)/W_1$ by mapping $\overline{G}_1^g$ to $W_1\pi(g^\sigma g^{-1})$.
\end{lemma}

\begin{lemma}\label{l:bijection2} Let $\overline G$ be a connected reductive
algebraic group and $\sigma$ be a Frobenius map of $\overline G$. Let $\overline G_1$ be a $\sigma$-stable reductive
subgroup of maximal rank of $\overline G$. Let $\overline G_1^g$ be also $\sigma$-stable. Denote $n=g^\sigma g^{-1}$.
Then $(\overline G_1^g)_\sigma=((\overline G_1)_{\sigma\circ n})^g$. Moreover, for arbitrary  $n\in
N_{\overline{G}}(\overline T)\cap N_{\overline{G}}(\overline G_1)$ group $(\overline G_1)_{\sigma\circ n}$ is conjugate
in $\overline G$ to some reductive subgroup of $\overline G_\sigma$.
\end{lemma}

\begin{prf} Let $h\in \overline G_1$. We have $h\in (\overline G_1)_{\sigma\circ n}$ if and only if
$h^{\sigma\circ n}=h$. By substituting the expression for $n$ we obtain $h^{\sigma\circ g^\sigma g^{-1}}=h$. The last
is equivalent to $h^{\sigma\circ g^\sigma}=h^g$. Whence it follows that $(h^g)^\sigma=h^g$, thus $h^g\in(\overline
G_1^g)_\sigma$. Therefore, $h$ lies in $(\overline G_1)_{\sigma\circ n}$ if and only if $h^g$ lies in $(\overline
G_1^g)_\sigma$. Hence $(\overline G_1^g)_\sigma=((\overline G_1)_{\sigma\circ n})^g$. By Lang --- Steinberg theorem
(see e.g. \cite[Theorem 10.1]{Steinberg}) for every element $n$ of $N_{\overline{G}}(\overline T)\cap
N_{\overline{G}}(\overline G_1)$ there exists an element $g$ of $\overline G$ such that $n=g^\sigma g^{-1}$. Thus the
second assertion of the lemma follows from the listed equivalences.
\end{prf}

Lemmas \ref{l:bijection} and \ref{l:bijection2} implies that a description of structure of reductive subgroups of a
finite group $\overline G_\sigma$ can be obtained by using the following scheme. We choose some set $\mathfrak{M}$
containing the full system of representatives of conjugacy classes of $\sigma$\nobreakdash-stable reductive subgroups
of maximal rank of $\overline G$. For every subgroup $\overline G_1$ in $\mathfrak{M}$ we describe the structure of
group of the form $(\overline G_1)_{\sigma\circ n}$, where $W_1\pi(n)$ runs through the full system of representatives
of conjugacy classes of $N_W(W_1)/W_1$. In what follows the set $\mathfrak{M}$ will be always a set of reductive
subgroups containing fixed $\sigma$-stable maximal torus of $\overline G$.

\begin{lemma}{\rm \cite[Proposition 0.5]{Tes}}\label{l:pperoid} Let $\overline G$ be a simple algebraic group
over an algebraically closed field of positive characteristic $p$ and let $\sigma$ be a Frobenius map of $\overline G$.
Then $p$-period of $\overline G_\sigma$ is equal to the minimal power of $p$ greater then the maximal height of a root
in root system of $\overline G$.
\end{lemma}

Remind that the maximal height of a root in systems of type $A_n$ is equal to $n-1$, in systems of types $B_n$ and
$C_n$ is equal to $2n-1$, and in system of type $D_n$ is equal to $2n-3$.

For a finite group $G$ denote by $\operatorname{exp}(G)$ the period of $G$. For a partition $\lambda=(\lambda_1,
\lambda_2,\dots,\lambda_k)$ denote by $|\lambda|$ the sum $\lambda_1+\lambda_2+\dots+\lambda_k$, the number $k$ is
called the length of $\lambda$. For a prime $p$ denote by $k_{\{p\}}$ the maximal power of $p$ dividing $k$. For
naturals $n_1$, $n_2,\dots, n_s$ denote $(n_1,n_2,\dots,n_s)$ their greatest common divisor and by
$[n_1,n_2,\dots,n_s]$ their least common multiple. Besides, the least common multiple of naturals $n_i$, where $i$ runs
through some set of indices, is denoted by $\underset{i}{\operatorname{lcm}}\{n_i\}$. Denote by
$\operatorname{diag}(A_1, A_2,\dots, A_k)$ a block-diagonal matrix with blocks $A_1$, $A_2, \dots, A_k$ and by
$\operatorname{diag}(a_1, a_2,\dots, a_k)$ a diagonal matrix with numbers $a_1$, $a_2, \dots, a_k$ on the diagonal.
Denote by $E_k$ the identity matrix of size $k\times k$.

We will also need the following lemma about abelian groups.

\begin{lemma}\label{l} Let $C_i$, $1\leqslant i\leqslant s$, be cyclic groups and
$B=C_1\times C_2\times\dots\times C_s$, where $s>1$. Let $A$ be a subgroup of $B$ of simple order $p$ such that $A\cap
C_i=1$ for all $1\leqslant i\leqslant s$. Then $\exp(B)=\exp(B/A)$.
\end{lemma}

\begin{prf} To prove the lemma it suffices to show that  $\exp(B)_{\{p\}}=\exp(B/A)_{\{p\}}$.
Let $j$ be such that $|C_j|_{\{p\}}\geqslant |C_i|_{\{p\}}$ for all $i$. We have $C_jA/A\simeq C_j/(C_j\cap A)\simeq
C_j$. The lemma is proved.
\end{prf}

The information about linear and unitary groups used in the paper can be found in \cite{Carter} or \cite{my}.

\sect{Spectra of symplectic groups}

In this section, $\overline G=Sp_{2n}(\overline F_p)$ is a symplectic group associated with form
$x_1y_{-1}-x_{-1}y_1+\ldots+x_ny_{-n}-x_{-n}y_n$, where $n\geqslant 1$. The map $\sigma$ acting on $\overline G$ as
follows: matrix $(a_{ij})$ under the action of $\sigma$ transforms into matrix $(a_{ij}^q)$, where $q$ is a power of
$p$, is a Frobenius map of $\overline G$. Moreover, $\overline{G}_\sigma=Sp_{2n}(q)$ and
$\overline{G}_\sigma/Z(\overline{G}_\sigma)=PSp_{2n}(q)$. Groups $Sp_2(q)$ and $PSp_2(q)$ are isomorphic to groups
$SL_2(q)$ and $PSL_2(q)$, description of their spectra can be found in \cite{my}.

If $p=2$ then, as it was mentioned above, $PSp_{2n}(q)\simeq\Omega_{2n+1}(q)$. In this case it is more convenient for
us to consider this groups as orthogonal, and a description of their spectra will be given in the following sections.
In this section we assume that characteristic $p$ is odd.

Rows and columns of a symplectic $2n$-matrix are numerated in the order $1,2,\ldots,n,-1,-2,\ldots,-n$. Denote by
$E_{i,j}$, where $i,j\in\{\pm 1,\pm2,\dots,\pm n\}$, a matrix containing identity at its $(i,j)$th place and containing
zeros at other places.

The maximal torus $\overline{T}$ of $\overline G$ consisting of all matrices of the form
$\operatorname{diag}(D,D^{-1})$, where $D$ is a nondegenerate diagonal matrix of size $n\times n$,
is $\sigma$-stable. Group $\overline N=N_{\overline{G}}(\overline T)$ is a subgroup of a group of
monomial matrices, thus, there exists a natural embedding of $W$ into the group of permutation of
the set $\{1,2,\ldots,n,-1,-2,\ldots,-n\}$. The image of the Weyl group $W\simeq\overline
N/\overline T$ under this embedding coincides with the group  $Sl_{n}$ of permutation $\tau$ such
that the equality $\tau(-i)=-\tau(i)$ holds. The map $\sigma$ act trivially on $W$, therefore,
$\sigma$\nobreakdash-conjugacy classes in $W$ coincide with usual conjugacy classes.

If we drop signs from elements of the set $\{1,2,\ldots,n,-1,-2,\ldots,-n\}$ we obtain a homomorphism from $Sl_n$ to
$Sym_n$. Assume $\tau\in Sl_n$ is mapped into a cycle $(i_1i_2\ldots i_k)$ then  $\tau$ fixes all elements distinct
from $\pm i_1, \pm i_2,\ldots, \pm i_k$. If $\tau^k(i_1)=i_1$ then $\tau$ is called {\it a positive cycle of length}
$k$; if $\tau^k(i_1)=-i_1$ then $\tau$ is called {\it a negative cycle of length} $k$. The image of an arbitrary
element $\tau$ of $Sl_n$ can be uniquely expressed as a product of disjoint cycles, and in accordance with this
factorization, $\tau$ can be uniquely expressed as a product of disjoint positive and negative cycles. Lengths of the
cycles together with their signs give a set of integers, called the cycle-type of $\tau$. Two elements of $Sl_n$ are
conjugate if they have the same cycle-type. We will denote the positive cycle of length $k$ by $(i_1 i_2\dots i_k)$ and
negative~--- by $\overline{(i_1 i_2\dots i_k)}$.

The root system $\Phi$ of $\overline G$ has type $C_n$ and can be represented as follows. Let $e_1,
\dots, e_n$ be a n orthonormal basis of a Euclidean space of dimension $n$. Then $\Phi=\{ \pm e_i
\pm e_j, \pm 2e_i, i\neq j, 1\leqslant i,j \leqslant n\}$. An element $w\in W$ acts on the basis
vectors in the following way: $w(e_i)=e_j$, if $w(i)=j$, and $w(e_i)=-e_j$, if $w(i)=-j$.

Define elements $u_r(t)$, where $r\in\Phi$ and $t\in\overline F_p$, as follows
$$u_r(t)=E_{2n}+t(E_{i,j}-E_{-j,-i}), \text{ if } r=e_i-e_j,$$
$$u_r(t)=E_{2n}+t(E_{i,-j}-E_{j,-i}), \text{ if } r=e_i+e_j,$$
$$u_r(t)=E_{2n}+t(E_{-i,j}-E_{-j,i}), \text{ if } r=-e_i-e_j,$$
$$u_r(t)=E_{2n}+t E_{i,-i}, \text{ if } r=2e_i,$$
$$u_r(t)=E_{2n}+t E_{-i,i}, \text{ if } r=-2e_i.$$

Then a root subgroup $\overline{U}_r$ ( with respect to torus $\overline T$) consists of elements $u_r(t)$,
$t\in\overline F_p$ ( see e.g. the proof of Theorem 11.3.2 from \cite{CarterBook2}). Remark that all root subgroups are
$\sigma$-stable. Thus all reductive subgroups containing $\overline T$ are $\sigma$-stable. In addition, an arbitrary
reductive subgroup of maximal rank is conjugate to some reductive subgroup containing $\overline T$. Therefore to
describe the structure of reductive subgroups of $\overline G_\sigma$ it suffices for every reductive subgroup
$\overline G_1$ containing $\overline T$, to describe a structure of a group of the form $(\overline G_1)_{\sigma\circ
n}$, where $W_1\pi(n)$ runs trough the full system of representatives of conjugacy classes of $N_W(W_1)/W_1$.

Subsystems $\Phi_1$ of $\Phi$ can be obtained in the following way (see \cite[section 3]{Carter}).
Let $\lambda$ and $\mu$ be partitions with $|\lambda|+|\mu|=n$. Let $\lambda=(\lambda_1, \lambda_2,
\dots, \lambda_l)$ and $\mu=(\mu_1,\mu_2,\dots, \mu_m)$. Let $I_1, I_2,\dots, I_l, J_1, J_2, \dots,
J_m$ be pairwise non-intersecting subsets of the set $\{1,2,\dots,n\}$ such that $|I_a|=\lambda_a$
for all $1\leqslant a\leqslant l$ and $|J_b|=\mu_b$ for all $1\leqslant b\leqslant m$. Define a
system $\Phi_1$ to be the union
$$\bigcup_{1\leqslant a\leqslant l}\{e_i-e_j : i\neq j, i,j\in
I_a\}\cup\bigcup_{1\leqslant b \leqslant m}\{\pm e_i \pm e_j, \pm 2e_i : i\neq j, i,j\in J_b\}.$$ Then $\Phi_1$ is a
root subsystem of $\Phi$ of type
$$A_{\lambda_1-1}\times A_{\lambda_2-1}\times \dots\times A_{\lambda_l-1}\times C_{\mu_1}\times
C_{\mu_2}\times \dots\times C_{\mu_m}$$ and every subsystem of $\Phi$ is equivalent to some system of such form under
the action of $W$.

For a subsystem $\Phi_1$ of $\Phi$, define an index set $I_{\Phi_1}$ using the following rule: indices $i$ and $-i$
belong to $I_{\Phi_1}$ if and only if $\Phi_1$ contains a root not orthogonal to $e_i$. We will denote by $\overline
T_{\Phi_1}$ a subgroup of $\overline T$ consisting of all diagonal matrices which have $1$'s at places with numbers not
in  $I_{\Phi_1}$ along the diagonal. Denote the group $\langle \overline T_{\Phi_1}, \overline U_r, r\in
{\Phi_1}\rangle$ by $\overline G_{\Phi_1}$.

Let ${\Phi_1}$ be an indecomposable subsystem and $k=|I_{\Phi_1}|/2$. If ${\Phi_1}$ contains only short roots, then it
has type $A_{k-1}$. By the choice of $\overline T$ and the form of root вида elements every matrix in $\overline
G_{\Phi_1}$ is block-diagonal with blocks $A$ and $A^{-\top}$ along the diagonal, where $A$ is a matrix of size $n
\times n$ and $-\top$ denotes a transpose-inverse map. After a suitable renumbering of rows and columns the matrix $A$
is also block-diagonal with blocks $E_{n-k}$ and $B$ along the diagonal, where $B$ is a matrix from $GL_k(\overline
F_p)$. If $A$ runs through all root elements of $\overline G_{\Phi_1}$, then $B$ runs through all root elements of
$GL_k(\overline F_p)$. In addition, if $A$ runs through $\overline T_{\Phi_1}$, then corresponding matrix $B$ runs
through the maximal torus of $GL_k(\overline F_p)$ consisting of all diagonal matrices. Since $GL_k(\overline F_p)$ is
generated by all its root subgroups and an arbitrary maximal torus, group $\overline G_{\Phi_1}$ is isomorphic to
$GL_k(\overline F_p)$. If ${\Phi_1}$ contains long roots, then it has type $C_k$. By applying arguments similar to
those above, one can show that $\overline G_{\Phi_1}$ is isomorphic to $Sp_{2k}(\overline F_p)$ is this case.

Let $\overline G_1$ be a $\sigma$-stable reductive subgroup of maximal rank with root system $\Phi_1$ of the form given
above. Denote by $n_{i}$ the number of sets of cardinality $i$ among the sets $I_a$ and by $m_j$ the number of sets of
cardinality $j$ among the sets $J_b$. Then $\overline G_1$ is isomorphic to direct product
$$\prod_i GL_i(\overline F_p)^{n_i}\times \prod_j Sp_{2j}(\overline
F_p)^{m_j},$$ and $N_W(W_1)/W_1$ is isomorphic to direct product
$$\prod_i  Sl_{n_i}\times \prod_j Sym_{m_j}$$(see \cite[Section 3]{Carter}).

Let us show that for every element $\overline w$ of $N_{W}(W_1)/W_1$ there exists a permutation matrix $\phi(\overline
w)$ such that groups, corresponding to the conjugacy class of $N_{W}(W_1)/W_1$ containing $\overline w$,  are conjugate
in $\overline G$ to group $(\overline G_1)_{\sigma\circ \phi(\overline w)}$. Denote by $\psi$ the natural homomorphism
from $N_W(W_1)$ to $N_W(W_1)/W_1$.

Let $\overline w$ be a positive cycle of length $k$ from direct factor $Sl_{n_i}$. Then $\overline G_1$ has a direct
factor isomorphic to a direct power of $GL_i(\overline F_p)^k$ consisting, after suitable renumbering of rows and
columns, of block-diagonal matrices of the form
$$\begin{pmatrix}X_1&&&&&&&&\\&\ddots&&&&&&&\\&&X_k&&&&&&\\&&&&E_{n-ki}&&&&\\&&&&&X_1^{-T}&&&\\&&&&&&\ddots&&
\\&&&&&&&X_k^{-T}&\\&&&&&&&&E_{n-ki}\end{pmatrix},$$ where $X_j\in GL_i(\overline F_p)$ for $1\leqslant j\leqslant k$.
The complete inverse $\overline w$ under $\psi$ contains a permutation $w$ equal to a product of positive cycles
$(1,i+1, 2i+1,\dots, (k-1)i+1)(2,i+2, 2i+2,\dots, (k-1)i+2)\dots(i,2i, 3i,\dots, ki)$. Let $\phi(\overline w)$ be the
permutation matrix corresponding to $w$. The matrix $\phi(\overline w)$ is symplectic and lies in $N_{\overline
G}(\overline G_1)\cap N_{\overline G}(\overline T)$. Let $\overline w$ be a negative cycle of length $k$. Then
$\overline G_1$ also contains a subgroup $GL_i(\overline F_p)^k$ of the form described above. The inverse image of
$\overline w$ under $\psi$ contains a permutation $w$ equal to a product of negative cycles $\overline{(1,i+1,
2i+1,\dots, (k-1)i+1)}\dots\overline{(i,2i, 3i,\dots, ki)}$. Let $\phi(\overline w)$ be the permutation matrix
corresponding to $w$. The matrix $\phi(\overline w)$ is not symplectic. Denote by $d$ the matrix
$\operatorname{diag}(E_{n+(k-1)i},-E_i, E_{n-ki})$. Then the matrix $d\phi(\overline w)$ is symplectic, lies in
$N_{\overline G}(\overline G_1)\cap N_{\overline G}(\overline T)$ and $\pi(d\phi(\overline w))=w$. Moreover, the action
of $\phi(\overline w)$ and $d\phi(\overline w)$ on $\overline G_1$ coincide. Thus, $(\overline G_1)_{\sigma\circ
d\phi(\overline w)}=(\overline G_1)_{\sigma\circ \phi(\overline w)}$. Let $\overline w$ be a cycle of length $k$ lying
in a direct factor $Sym_{m_i}$ of $N_W(W_1)/W_1$. Then $\overline G_1$ contains a subgroup isomorphic to
$Sp_{2i}(\overline F_p)^k$, consisting, after suitable renumbering of rows and columns, of block-diagonal matrices of
the form $\operatorname{diag}(X_1, X_2,\dots,X_k, E_{2n-2ki})$, $X_j\in Sp_{2i}(\overline F_p)$ for $1\leqslant
j\leqslant k$. The matrix $\phi(\overline w)$ is defined in the same way as in the positive cycle case. By extending
$\phi$ to $N_{W}(W_1)/W_1$, we obtain the desirable map.

The conjugacy classes of $N_W(W_1)/W_1$ can be determined by collections of partitions. Indeed, a
conjugacy class of $Sym_{m_j}$ is determined by a partition $\rho^{(j)}$ of $m_j$; a conjugacy
class of $Sl_{n_i}$ is determined by a pair of partitions $\xi^{(i)}$ and $\zeta^{(i)}$ with
$|\xi^{(i)}|+|\zeta^{(i)}|=n_i$, where parts of $\xi^{(i)}$ give lengths of positive cycles and
parts of $\zeta^{(i)}$ give lengths of negative cycles. We have
\begin{equation}\label{equation}\sum_i i(|\xi^{(i)}|+|\zeta^{(i)}|)+\sum_j j|\rho^{(j)}|=n.
\end{equation}

We carry out the description of reductive subgroups of maximal rank of finite group $\overline G_\sigma$ according to
the following scheme. We shall choose a collection of partitions $\rho^{(i)}, \xi^{(i)}, \zeta^{(i)}$ satisfying
(\ref{equation}). It determines up to equivalence under the action of $W$ the root system of $\overline G_1$: the sum
$|\xi^{(i)}|+|\zeta^{(i)}|$ determines the number of subsystems of type $A_{i-1}$ and number $|\rho^{(i)}|$ determines
the number of subsystems of type $C_i$. It also determines the conjugacy class of $N_W(W_1)/W_1$. We shall choose a
representative of this class $\overline w$ and describe a structure of group $(\overline
G_1)_{\sigma\circ\phi(\overline w)}$.

The element $\overline w$ can be presented as a product of disjoint cycles, and $\overline G_1$ can
be presented as a direct product of subgroups such that the image of each cycle under $\phi$ acts
nontrivially exactly on one factor. Thus, to describe the structure of reductive subgroups it
suffices, for every cycle, to describe the structure of corresponding subgroup. Group
$N_W(W_1)/W_1$ contains three types of cycles: positive and negative cycles from factors $Sl_{n_i}$
and cycles from factors $Sym_{m_i}$.

Let $\xi^{(i)}_j$ be a part of partition $\xi^{(i)}$. Then the decomposition of $\overline w$ into disjoint cycles
contains a positive cycle of length $\xi^{(i)}_j$. This implies that $\overline G_1$ contains a subgroup isomorphic to
$GL_i(\overline F_p)^{\xi^{(i)}_j}$ consisting of block-diagonal matrices of the form
$$\operatorname{diag}(X_1,X_2 ,\dots, X_{\xi_j^{(i)}}, E, X_1^{-\top},X_2^{-\top} ,\dots, X_{\xi_j^{(i)}}^{-\top}, E),$$
where $X_k\in GL_i(\overline F_p)$ for $1\leqslant k\leqslant \xi^{(i)}_j$ and $E$ is an identity matrix of suitable
size. We have
\begin{equation*}\begin{split}&\sigma\circ\phi(\overline w)(\operatorname{diag}(X_1,X_2 ,\dots, X_{\xi_j^{(i)}}, E,
X_1^{-\top},X_2^{-\top} ,\dots, X_{\xi_j^{(i)}}^{-\top}, E))=\\&=\operatorname{diag}(X_1,X_2 ,\dots, X_{\xi_j^{(i)}},
E, X_1^{-\top},X_2^{-\top} ,\dots, X_{\xi_j^{(i)}}^{-\top}, E).\end{split}\end{equation*}

Rewrite \begin{equation*}
\begin{split}
&\sigma(\operatorname{diag}(X_{\xi_j^{(i)}},X_1 ,\dots, X_{\xi_j^{(i)}-1}, E, X_{\xi_j^{(i)}}^{-\top},X_1^{-\top} ,\dots,
X_{\xi_j^{(i)}-1}^{-\top}, E))=\\&=\operatorname{diag}(X_1,X_2 ,\dots, X_{\xi_j^{(i)}}, E, X_1^{-\top},X_2^{-\top} ,\dots,
X_{\xi_j^{(i)}}^{-\top}, E),
\end{split}
\end{equation*}
this yields the system of equalities $X_2=X_1^\sigma$, $X_3=X_2^\sigma, \dots,
X_1=X_{\xi_j^{(i)}}^\sigma$. This system gives equation $X_1=X_1^{\sigma^{\xi^{(i)}_j}}$, which is
true if and only if $X_1$ lies in $GL_{i}(q^{\xi^{(i)}_j})$. Thus, in matrices of group $(\overline
G_1)_{\sigma\circ \phi(\overline w)}$, parts $\xi^{(i)}_j$ correspond to blocks of the form
$\operatorname{diag}(Y, Y^{-\top})$, where
$$Y=\operatorname{diag}(X, X^\sigma,\dots, X^{\sigma^{(\xi^{(i)}_j-1)}}),\; X\in GL_{i}(q^{\xi^{(i)}_j}).$$ Denote
the group of matrices of the form $\operatorname{diag}(Y, Y^{-\top})$ by $H_{ij}^\xi$.

Let $\zeta^{(i)}_j$ be a part of partition $\zeta^{(i)}$. Then the decomposition of $\overline w$ into disjoint cycles
contains a negative cycle of length $\zeta^{(i)}_j$. This implies that $\overline G_1$ contains a subgroup isomorphic
to $GL_i(\overline F_p)^{\zeta^{(i)}_j}$ consisting of block-diagonal matrices of the form
$$\operatorname{diag}(X_1,X_2 ,\dots, X_{\zeta_j^{(i)}}, E, X_1^{-\top},X_2^{-\top} ,\dots, X_{\zeta_j^{(i)}}^{-\top}, E),$$
where $X_k\in GL_i(\overline F_p)$ for $1\leqslant k\leqslant \zeta^{(i)}_j$ and $E$ is an identity matrix of suitable
size. We have \begin{equation*}\begin{split}&\sigma\circ\phi(\overline w)(\operatorname{diag}(X_1,X_2 ,\dots,
X_{\zeta_j^{(i)}}, E, X_1^{-\top},X_2^{-\top} ,\dots, X_{\zeta_j^{(i)}}^{-\top}, E))=\\&=\operatorname{diag}(X_1,X_2
,\dots, X_{\zeta_j^{(i)}}, E, X_1^{-\top},X_2^{-\top} ,\dots, X_{\zeta_j^{(i)}}^{-\top}, E).\end{split}\end{equation*}
Rewrite \begin{equation*}
\begin{split}
&\sigma(\operatorname{diag}(X_{\zeta_j^{(i)}}^{-\top},X_1 ,\dots, X_{\zeta_j^{(i)}-1}, E, X_{\zeta_j^{(i)}},X_1^{-\top} ,\dots,
X_{\zeta_j^{(i)}-1}^{-\top}, E))=\\&=\operatorname{diag}(X_1,X_2 ,\dots, X_{\zeta_j^{(i)}}, E, X_1^{-\top},X_2^{-\top} ,\dots,
X_{\zeta_j^{(i)}}^{-\top}, E),
\end{split}
\end{equation*}
this yields the system of equalities $X_2=X_1^\sigma$, $X_3=X_2^\sigma, \dots, X_1=(X_{\zeta_j^{(i)}}^{-\top})^\sigma$.
This system gives an equation $X_1=(X_1^{-\top})^{\sigma^{\zeta^{(i)}_j}}$, which is true if and only if $X_1$ lies in
$GU_{i}(q^{\zeta^{(i)}_j})$. Thus, in matrices of group $(\overline G_1)_{\sigma\circ \phi(\overline w)}$ parts
$\zeta^{(i)}_j$ correspond to blocks of the form $\operatorname{diag}(Y, Y^{-\top})$, where
$$Y=\operatorname{diag}(X, X^\sigma,\dots, X^{\sigma^{(\zeta^{(i)}_j-1)}}),\; X\in GU_{i}(q^{\zeta^{(i)}_j}).$$ Denote
the group of matrices of the form $\operatorname{diag}(Y, Y^{-\top})$ by $H_{ij}^\zeta$.

By similar arguments we can show that parts $\rho^{(i)}_j$ of partition $\rho^{(i)}$ correspond to blocks of the form
$$Y=\operatorname{diag}(X, X^\sigma, X^{\sigma^2},\dots, X^{\sigma^{(\rho^{(i)}_j-1)}}), \;X\in
Sp_{2i}(q^{\rho^{(i)}_j}).$$ Denote the group of matrices of such form by $H_{ij}^\rho$.

Observe that groups $H_{1j}^\xi$ and $H_{1j}^\zeta$ consist of diagonal matrices. Moreover, a group of diagonal
matrices with blocks $A_j$ and $B_j$ along the diagonal, where $A_j\in H_{1j}^\xi$ and $B_j\in H_{1j}^\zeta$, in
contained in $Sp_{2n_1}(\overline F_p)$ and conjugate in it to a maximal torus of $Sp_{2n_1}(q)$, corresponding to the
pair of partitions $\xi^{(1)}$ and $\zeta^{(1)}$ (see \cite[Proposition 3.1]{our}). Thus, we proved

\begin{prop}\label{p:spfin} Let $\overline G_1^g$ be a $\sigma$-stable subgroup of $\overline G$. Let $w=\pi(g^\sigma
g^{-1})$ and coset $W_1w$ lie in the conjugacy class corresponding to a collection of partitions $\rho^{(i)},
\xi^{(i)}, \zeta^{(i)}$. Then $(\overline G_1^g)_\sigma$ is isomorphic to group $G_1$ consisting of all block-diagonal
matrices with blocks $A_{ij}$, $B_{ij}$, $C_{ij}$ along the diagonal, where $A_{ij}\in H_{ij}^\xi$, $B_{ij}\in
H_{ij}^\zeta$ and $C_{ij}\in H_{ij}^\rho$. In particular,
$$\overline G_1^g\simeq
\prod_{i,j}GL_i(q^{\xi^{(i)}_j})\times\prod_{i,j}GU_i(q^{\zeta^{(i)}_j})\times\prod_{i,j}Sp_{2i}(q^{\rho^{(i)}_j}).$$
\end{prop}

Since the center of $GL_{i}(q^{\xi^{(i)}_j})$ is a cyclic group of order $q^{\xi^{(i)}_j}-1$, the center of
$GU_{i}(q^{\zeta^{(i)}_j})$ ia a cyclic group of order $q^{\zeta^{(i)}_j}+1$, the center of group $G_1$ from the
statement of Proposition \ref{p:spfin} is
\begin{equation}\label{zsp}\prod_{i,j}(q^{\xi^{(i)}_j}-1)\times\prod_{i,j}(q^{\zeta^{(i)}_j}+1)
\times\prod_{i,j}Z(Sp_{2i}(q^{\rho^{(i)}_j})).\end{equation}

The following assertion gives a description of composite part of spectrum of $Sp_{2n}(q)$.

\begin{theorem}\label{t:sp} Let $G\simeq Sp_{2n}(q)$, where $n\geqslant 2$ and $q$ is a power of add prime $p$.
Assume that for every natural $k$ such that $2n_0=p^{k-1}+1< 2n$ and for every pair of partitions $\alpha=(\alpha_1,
\alpha_2,\dots, \alpha_a)$ and $\beta=(\beta_1, \beta_2,\dots, \beta_b)$ such that $n-n_0=|\alpha|+|\beta|$, the set
$\nu(G)$ contains a number
$$p^k[q^{\alpha_1}-1, q^{\alpha_2}-1, \dots, q^{\alpha_a}-1, q^{\beta_1}+1, q^{\beta_2}+1,\dots, q^{\beta_b}+1],$$
contains $2p^k$ if $2n=p^{k-1}+1$ for some $k>1$ and does not contain any other number. Then $\mu_m(G)\subseteq \nu(G)
\subseteq \omega(G)$.
\end{theorem}

\begin{prf}
We shall prove the inclusion $\nu(G)\subseteq \omega(G)$ first. To do this we for every possible
$k$ and partitions $\alpha$ and $\beta$ shall construct a reductive subgroup $H$ of $G$ such that
$\eta(H)$ is equal to corresponding element of $\nu(G)$. Let $I=\{1, 2,\dots, n_0\}$. Put
$\Phi_1=\{\pm e_i \pm e_j, \pm 2e_i, i\neq j, i,j\in I\}$. Then $\Phi_1$ is a subsystem of $\Phi$
of type $C_{n_0}$. Define $\overline G_1=\langle \overline T, \overline U_r, r\in\Phi_1\rangle$.
Let $W_1$ a the Weyl group of $\overline G_1$. Group $N_W(W_1)/W_1$ is isomorphic to $Sl_{n-n_0}$.
Consider the conjugacy class of $N_W(W_1)/W_1$ corresponding to a pair of partitions $\alpha$ and
$\beta$, where parts of $\alpha_i$ give lengths of positive cycles and parts of $\beta_j$ give
lengths of negative cycles. Group $G$ contains a reductive subgroup of maximal rank $(\overline
G_1^g)_\sigma$ which structure is determined by this class. Let $a$ be the number of elements in
$\alpha$ and $b$ be the number of elements in $\beta$. By Proposition \ref{p:spfin}
$$(\overline G_1^g)_\sigma\simeq Sp_{2n_0}(q)\times(q^{\alpha_1}-1)\times \dots\times
(q^{\alpha_a}-1)\times(q^{\beta_1}+1)\times\dots\times(q^{\beta_b}+1).$$ By Lemma \ref{l:pperoid} the maximal order of
unipotent element in $Sp_{2n_0}(q)$ is equal to $p^k$. Therefore, $\eta((\overline G_1^g)_\sigma)=p^k[q^{\alpha_1}-1,
q^{\alpha_2}-1, \dots, q^{\alpha_a}-1, q^{\beta_1}+1, q^{\beta_2}+1,\dots, q^{\beta_b}+1]$.

Moreover, if $2n=p^{k-1}+1$, then $p^k$ is the maximal power of $p$ lying in $\omega(G)$, and $p^k$ is not contained in
spectrum of any proper reductive subgroup. Thus, the number $\eta(G)$ equal to $2p^k$ lies in $\mu_m(G)$ and,
therefore, must lie in $\nu(G)$. If $2n\neq p^{k-1}+1$ for every $k\geqslant 1$, then $\eta(G)$ does not lie in
$\mu(G)$ and divides $\eta(H)$ for some proper reductive subgroup of $G$.

It remains to prove the inclusion $\mu_m(G)\subseteq \nu(G)$. To do this we shall show that for
every proper reductive subgroup of maximal rank $H$, which is not a maximal torus, there exists an
element of $\nu(G)$ divided by $\eta(H)$. Let the structure of $H$ be determined by a collection of
partitions $\xi^{(i)}, \zeta^{(i)}$, $1\leqslant i\leqslant a$ and $\rho^{(i)}$, $1\leqslant i
\leqslant b$. Since $H$ is not a maximal torus, one of inequalities $a>1$ or $b>0$ hods. By formula
(\ref{zsp}) the center of $H$ is isomorphic to the direct product
$$\prod_{i,j}(q^{\xi^{(i)}_j}-1)\times\prod_{i,j}(q^{\zeta^{(i)}_j}+1)\times\prod_{i,j}Z(Sp_{2i}(q^{\rho^{(i)}_j})).$$
Put $e=\exp(\prod_{i,j}Z(Sp_{2i}(q^{\rho^{(i)}_j})))$. We have
$$\eta(H)=p^k[\underset{i,j}{\operatorname{lcm}}\{q^{\xi^{(i)}_j}-1\}, \underset{i,j}{\operatorname{lcm}}
\{q^{\zeta^{(i)}_j}+1\}, e]$$ here $k$ is such that $p^k$ is the minimal power of $p$ greater then
$\operatorname{max}\{a-1, 2b-1\}$, i.\,e. the greatest power of $p$ lying n $\omega(H)$. Let
$2n_0=p^{k-1}+1$. Then
$$n_0=\frac{p^{k-1}+1}{2}\leqslant\frac{\max\{a-1, 2b-1\}+1}{2}=\max\{a/2,b\}.$$

By equality (\ref{equation}) we have
$$|\xi^{(1)}|+|\zeta^{(1)}|+2(|\xi^{(2)}|+|\zeta^{(2)}|)+\dots+a(|\xi^{(a)}|+|\zeta^{(a)}|)+$$
$$+|\rho^{(1)}|+2|\rho^{(2)}|+\dots+b|\rho^{(b)}|=n.$$
Put
$$x=|\xi^{(1)}|+|\zeta^{(1)}|+|\xi^{(2)}|+|\zeta^{(2)}|+\dots+|\xi^{(a)}|+|\zeta^{(a)}|+n_0.$$ Let us show that
$x\leqslant n$. By substituting the expressions for $x$ and $n$, we obtain the inequality
$$|\xi^{(1)}|+|\zeta^{(1)}|+|\xi^{(2)}|+|\zeta^{(2)}|+\dots+|\xi^{(a)}|+|\zeta^{(a)}|+n_0\leqslant$$
$$\leqslant|\xi^{(1)}|+|\zeta^{(1)}|+2(|\xi^{(2)}|+|\zeta^{(2)}|)+\dots+a(|\xi^{(a)}|+|\zeta^{(a)}|)+$$
$$+|\rho^{(1)}|+2|\rho^{(2)}|+\dots+b|\rho^{(b)}|.$$ Rewrite $$n_0\leqslant
|\xi^{(2)}|+|\zeta^{(2)}|+\dots+(a-1)(|\xi^{(a)}|+|\zeta^{(a)}|)+$$
$$+|\rho^{(1)}|+2|\rho^{(2)}|+\dots+b|\rho^{(b)}|.$$
If $\max\{a/2,b\}=b$, then $n_0\leqslant b|\rho^{(b)}|$ and, therefore, $x\leqslant n$. Assume that
$\max\{a/2,b\}=a/2$. Then $a\geqslant 2$ and required inequality is follows from inequalities
$n_0\leqslant a/2\leqslant (a-1)(|\xi^{(a)}|+|\zeta^{(a)}|).$

Define partitions $\alpha$ and $\beta$ as follows:
$$\alpha=(\xi^{(1)}_1, \xi^{(1)}_2,\dots, \xi^{(2)}_1, \xi^{(2)}_2,\dots, \xi^{(a)}_1, \dots, n-x),$$
$$\beta=(\zeta^{(1)}_1, \zeta^{(1)}_2,\dots, \zeta^{(2)}_1, \zeta^{(2)}_2,\dots, \zeta^{(a)}_1, \dots).$$

Let $\Phi_1$ be a subsystem of $\Phi$ of type $C_{n_0}$. Group $N_W(W_1)/W_1$ is isomorphic to
$Sl_{n-n_0}$. Choose in $N_W(W_1)/W_1$ the conjugacy class corresponding to the pair of partitions
$\alpha$ and $\beta$, where parts of $\alpha$ give the lengths of positive cycles and parts of
$\beta$ give the lengths of negative cycles. By Proposition \ref{p:spfin} a reductive subgroup
$G_1$ with the root system $\Phi_1$ corresponding to this class is isomorphic to
$$Sp_{2n_0}(q)\times\prod_{i,j}(q^{\xi^{(i)}_j}-1)\times\prod_{i,j}(q^{\zeta^{(i)}_j}+1)\times(q^{n-x}-1).$$
Hence
$$\eta(G_1)=p^k[2,\underset{i,j}{\operatorname{lcm}}\{q^{\xi^{(i)}_j}-1\},
\underset{i,j}{\operatorname{lcm}} \{q^{\zeta^{(i)}_j}+1\}, q^{n-x}-1].$$ This number is divided by
$\eta(H)$ and lies in $\nu(G)$. The theorem is proved.
\end{prf}

\begin{theorem}\label{t:psp} Let $G\simeq PSp_{2n}(q)$, where $n\geqslant 2$ and $q$ is a power of odd prime $p$.
Assume that for every natural $k$ such that $2n_0=p^{k-1}+1< 2n$ and for every pair of partritions
$\alpha=(\alpha_1, \alpha_2,\dots, \alpha_a)$ and $\beta=(\beta_1, \beta_2,\dots, \beta_b)$ such
that $n-n_0=|\alpha|+|\beta|$ the set $\nu(G)$ contains a number
$$p^k[q^{\alpha_1}-1, q^{\alpha_2}-1, \dots, q^{\alpha_a}-1, q^{\beta_1}+1, q^{\beta_2}+1,\dots, q^{\beta_b}+1],$$
and contains no other numbers. Then $\mu_m(G)\subseteq \nu(G) \subseteq \omega(G)$.
\end{theorem}

\begin{prf} The proof of Theorem \ref{t:sp} implies that the set $\nu(Sp_{2n}(q))$ consists of numbers $\eta(H)$,
where $H$ runs through a set of reductive subgroups of the form
$$Sp_{2n_0}(q)\times(q^{\alpha_1}-1)\times (q^{\alpha_2}-1)\times\dots\times
(q^{\alpha_a}-1)\times(q^{\beta_1}+1)\times(q^{\beta_2}+1)\times\dots\times(q^{\beta_b}+1),$$ where
$2n_0=p^{k-1}+1<2n$. If $H\neq Sp_{2n}(q)$, then $\eta(H/Z(Sp_{2n}(q)))=\eta(H)$  by Lemma \ref{l}.
If $H=Sp_{2n}(q)$, then $2n=p^{k-1}+1$ for some $k$ and $\eta(H/Z(Sp_{2n}(q))=p^k$. Therefore, sets
$\mu_m(Sp_{2n}(q))$ and $\mu_m(PSp_{2n}(q))$ are distinct if and only if $2n=p^{k-1}+1$ for some
$k\geqslant 1$. In the latter case $\mu_m(Sp_{2n}(q))=\mu_m(PSp_{2n}(q))\cup\{2p^k\}$.
\end{prf}

As a corallary of Lemma \ref{l:pperoid}, \cite[Proposition 3.1, Theorem 3]{our}, Theorems
\ref{t:sp} and \ref{t:psp} we obtain a description of spectra of finite symplectic and projective
symplectic  groups over field of odd characteristic.

\begin{cor}\label{c:specsp}  Let $G=Sp_{2n}(q)$, where $n\geqslant 2$ and $q$ is a power of odd prime $p$.
Then $\omega(G)$ consist of all divisor of the following numbers:
\begin{itemize}

\item[$1)$] $[q^{n_1}+\varepsilon_11,q^{n_2}+\varepsilon_21,\dots,q^{n_s}+\varepsilon_s1]$ for every $s\geqslant 1$,
$\varepsilon_i\in\{+,-\}$, $1\leqslant i\leqslant s$, and $n_1,n_2,\dots,n_s>0$ such that
$n_1+n_2+\dots+n_s=n$;

\item[$2)$] $p^k[q^{n_1}+\varepsilon_11,q^{n_2}+\varepsilon_21,\dots,q^{n_s}+\varepsilon_s1]$ for
every $s\geqslant 1$, $\varepsilon_i\in\{+,-\}$, $1\leqslant i\leqslant s$, and
$k,n_1,n_2,\dots,n_s>0$ such that $p^{k-1}+1+2n_1+2n_2+\dots+2n_s=2n$;

\item[$3)$] $2p^k$, if $p^{k-1}+1=2n$ for some $k>1$.
\end{itemize}
\end{cor}

\begin{cor}\label{c:specpsp}  Let $G=PSp_{2n}(q)$, where $n\geqslant 2$ and $q$ is a power of odd prime $p$.
Then $\omega(G)$ consist of all divisor of the following numbers:
\begin{itemize}
\item[$1)$] $\frac{q^n\pm1}{2}$;

\item[$2)$] $[q^{n_1}+\varepsilon_11,q^{n_2}+\varepsilon_21,\dots,q^{n_s}+\varepsilon_s1]$ for every $s\geqslant 2$,
$\varepsilon_i\in\{+,-\}$, $1\leqslant i\leqslant s$, and $n_1,n_2,\dots,n_s>0$ such that
$n_1+n_2+\dots+n_s=n$;

\item[$3)$] $p^k[q^{n_1}+\varepsilon_11,q^{n_2}+\varepsilon_21,\dots,q^{n_s}+\varepsilon_s1]$ for
every $s\geqslant 1$, $\varepsilon_i\in\{+,-\}$, $1\leqslant i\leqslant s$, and
$k,n_1,n_2,\dots,n_s>0$ such that $p^{k-1}+1+2n_1+2n_2+\dots+2n_s=2n$;

\item[$4)$] $p^k$, if $p^{k-1}+1=2n$ for some $k>1$.
\end{itemize}
\end{cor}

\sect{Reductive subgroups of orthogonal groups}

Recall some definitions and notations concerning orthogonal groups. Group $GO_{n}(K,Q)$ is an
orthogonal group of dimension $n$ over a field $K$, associated to a non-singular quadratic form
$Q$, and group $SO_n(K,Q)$ is a subgroup of $GO_n(K,Q)$ consisting of all matrices of determinant
$1$. If $K$ is algebraically closed, then all forms lead to isomorphic groups.

Let $K=GF(q)$, where $q$ is a power of a prime $p$. In the case of odd dimension all quadratic
forms lead to the same orthogonal group $GO_{n}(q)$; in the case of even dimension there are two
non-isomorphic groups $GO^+_{n}(q)$ and $GO^-_{n}(q)$. In all cases we will use unified notation
$GO^\varepsilon_n(q)$, where $\varepsilon$ is an empty symbol if $n$ is odd and
$\varepsilon\in\{+,-\}$ if $n$ is even.

Denote by $\Omega_n^\varepsilon(q)$ the commutator subgroup of $SO^{\varepsilon}_{n}(q)$. Group
$\Omega_{2n+1}(q)$ is simple if pair $(n,q)$ is distinct from $(1,2),(1,3),(2,2)$. Group
$\Omega_{2n}^\varepsilon(q)$ has the center of order $(4,q^{n}-\nolinebreak\varepsilon1)/2$ and its
factor group by the center is denoted by $P\Omega_{2n}^\varepsilon(q)$. For $n<4$ simple groups
$P\Omega_{2n}^\varepsilon(q)$ are isomorphic to some linear or unitary groups. For $n\geqslant 4$
all groups $P\Omega_{2n}^\varepsilon(q)$ are simple.

Put $\overline G=SO_{2n+1}(\overline F_p,Q)$, where $Q(x)=x_0^2+x_1x_{-1}+\dots+x_nx_{-n}$,
$n\geqslant 1$. Rows and columns of matrices from $\overline G$ are numerated in the order $0,
1,2,\ldots,n,-1,-2,\ldots,-n$. Define group $\overline H$ as a subgroup of $\overline G$ consisting
of all matrices of the form $\operatorname{diag}(1,A)$, where $A$ is a matrix of size $2n\times
2n$. Then $\overline H\simeq SO_{2n}(\overline F_p)$. Let map $\sigma$ act on $\overline G$ as
follows: matrix $(a_{ij})$ under action of $\sigma$ transforms into matrix $(a_{ij}^q)$, where $q$
is a power of $p$. Let $\tau$ be a graph automorphism of $\overline H$. Then $\sigma$ and
$\sigma\circ\tau$ are Frobenius maps of $\overline G$. In this notation $\overline G_\sigma\simeq
SO_{2n+1}(q)$, $\overline H_\sigma\simeq SO_{2n}^+(q)$ and $\overline H_{\sigma\circ\tau}\simeq
SO_{2n}^-(q)$.

The subgroup $\overline T$ of $\overline G$ consisting of all diagonal matrices of the form
$\operatorname{diag}(1,D,D^{-1})$ is a maximal torus of groups $\overline G$ and $\overline H$. The
Weyl group $W = N_{\overline G}(\overline T)/\overline T$ of $\overline G$ is isomorphic to $Sl_n$.
Its subgroup $W_{\overline H}=N_{\overline H}(\overline T)/\overline T$ is isomorphic to the
subgroup of $Sl_n$ consisting of permutations whose decomposition into disjoint cycles contain even
number of negative cycles. Put $n_0=\operatorname{diag}(-1,A_0)$, where $A_0$ is the permutation
matrix corresponding to negative cycle $(n,-n)$. Then $n_0\in N_{\overline G}(\overline T)$ and
$W=W_{\overline H} \cup w_0W_{\overline H}$, where $w_0=\pi(n_0)$.

The map $\sigma$ act trivially on $W$. Hence, $\sigma$-conjugacy classes of $W$ coincide with
conjugacy classes. Consider $\sigma\circ\tau$-conjugacy classes of $W_{\overline H}$. The action of
$\tau$ on $W_{\overline H}$ coincide with the action of $w_0$. Elements $w_1$ and $w_2$ of
$W_{\overline H}$ are $\sigma\circ\tau$-conjugate in $W_{\overline H}$ if
$w_1=(w^{-1})^{\sigma\circ\tau}w_2w$ for some $w$ of $W_{\overline H}$. Since
$w^{\sigma\circ\tau}=w_0^{-1}ww_0$, the equality $w_1= (w^{-1})^{\sigma\circ\tau}w_2w$ is
equivalent to $w_0w_1=w^{-1}w_0w_2w$. Thus $w_1$ and $w_2$ are $\sigma\circ\tau$-conjugate in
$W_{\overline H}$ if and only if $w_0w_1$ and $w_0w_2$ are conjugated by an element of
$W_{\overline H}$. For convenience, we will say that the structure of reductive subgroup of
$\overline H_{\sigma\circ\tau}$ is determined by a conjugacy class of $W$ contained in
$w_0W_{\overline H}$.

The root system $\Phi$ of $\overline G$ has type $B_n$. The system $\Phi$ can be presented as
follows. Let $e_1, e_2, \dots, e_n$ be an orthonormal basis of a Euclidian space of dimension $n$.
Then $\Phi=\{\pm e_i\pm e_j, \pm e_j, i\neq j, 1\leqslant i,j\leqslant n\}$. The group $W$ acts on
the basis vectors in the following way: $w(e_i)=e_j$ if $w(i)=j$ and $w(e_i)=-e_j$ if $w(i)=-j$.

The root system $\Phi_{\overline H}$ of $\overline H$ has type $D_n$. The system $\Phi_{\overline
H}$ can be presented as follows. Let $e_1, e_2, \dots, e_n$ be an orthonormal basis of a Euclidian
space of dimension $n$. Then $\Phi_{\overline H}=\{\pm e_i\pm e_j, i\neq j, 1\leqslant i,j\leqslant
n\}$. The group $W$ acts on the basis vectors in the following way: $w(e_i)=e_j$ if $w(i)=j$ and
$w(e_i)=-e_j$ if $w(i)=-j$. Thus the system $\Phi_{\overline H}$ is a subsystem of $\Phi$.

Let $p$ be odd. Define elements $u_r(t)$, $r\in \Phi,$ $t\in \overline{F}_p$ as follows:

$$u_r(t)=E+t(E_{ij}-E_{-j,-i}), \text{ if } r=e_i-e_j,$$
$$u_r(t)=E+t(E_{i,-j}-E_{j,-i}), \text{ if } r=e_i+e_j,$$
$$u_r(t)=E+t(E_{-j,i}-E_{-i,j}), \text{ if } r=-e_i-e_j,$$
$$u_r(t)=E+t(2E_{i,0}-E_{0,-i})-t^2E_{i,-i}, \text{ if } r=e_i,$$
$$u_r(t)=E-t(2E_{-i,0}-E_{0,i})-t^2E_{-i,i}, \text{ if } r=-e_i.$$

Let now $p=2$. In this case elements $u_r(t)$ are defined in the following way:

$$u_r(t)=E+t(E_{ij}-E_{-j,-i}), \text{ if } r=e_i-e_j,$$
$$u_r(t)=E+t(E_{i,-j}-E_{j,-i}), \text{ if } r=e_i+e_j,$$
$$u_r(t)=E+t(E_{-j,i}-E_{-i,j}), \text{ if } r=-e_i-e_j,$$
$$u_r(t)=E+t^2E_{i,-i}, \text{ if } r=e_i,$$
$$u_r(t)=E+t^2E_{-i,i}, \text{ if } r=-e_i.$$

In both cases root subgroup $\overline U_r$ of $\overline G$ and $\overline H$ (with respect to
$\overline T$) consist of elements $u_r(t)$, where $t\in \overline{F}_p$ ( see, e.g., the proof of
Theorem 11.3.2 from \cite{CarterBook2}). Here, in the case of group $\overline G$ roots $r$ lie in
the system $\Phi$ and in the case of group $\overline H$~--- in its subsytem $\Phi_{\overline H}$.
The root subgroups are $\sigma$-stable and, therefore, every reductive subgroup containing
$\overline T$ is $\sigma$-stable. The element $w_0$ defined above transforms the root $e_{n-1}-e_n$
into the root $e_{n-1}+e_n$ and fixes the rest of fundamental roots. Thus the action of $n_0$ on
$\overline H$ coincide with the action of graph automorphism.

The following isomorphisms $$\langle U_r, r\in\Phi\rangle\simeq\Omega_{2n+1}(\overline F_p),$$
$$\langle U_r, r\in\Phi_{\overline H}\rangle\simeq\Omega_{2n}(\overline F_p)$$ holds
(see \cite[Section 11.3]{CarterBook2}).

If $p$ is odd, then subsystems of $\Phi$ ca be obtained in the following way (see \cite[Section
3]{Carter}). Let $\lambda$ and $\mu$ are partitions satisfying the condition $|\lambda|+|\mu|\leq
n$ with no part of $\mu$ equal to $1$. Let $\lambda=(\lambda_1,\lambda_2,\dots,\lambda_l)$ and
$\mu=(\mu_1,\mu_2,\dots,\mu_m)$. Put $\rho=n-|\lambda|-|\mu|$. Let $I_1,I_2,\dots, I_l,
J_1,J_2,\dots, J_m$ be pairwise non-intersecting subsets of the set $\{1,2,\dots,n\}$ such that
$|I_a|=\lambda_a$ for all $1\leqslant a\leqslant l$ and $|J_b|=\mu_b$ for all $1\leqslant
b\leqslant m$, and let $K$ be the supplement to the union of the sets $I_1,I_2,\dots,I_l,
J_1,J_2,\dots, J_m$ in $\{1,2,\dots,n\}$. Define the system $\Phi_1$ to be the union
$$\bigcup_{1\leqslant a\leqslant l}\{e_i-e_j: i\neq j, i,j\in I_a\}\cup\bigcup_{1\leqslant b\leqslant m}\{\pm
e_i\pm e_j: i\neq j, i,j\in J_b\}\cup$$$$\cup\{\pm e_i\pm e_j, \pm e_i: i\neq j, i,j\in K\}.$$ Then
$\Phi_1$ is a subsystem of $\Phi$ of type $$A_{\lambda_1-1}\times A_{\lambda_2-1}\times\dots\times
A_{\lambda_l-1}\times D_{\mu_1}\times D_{\mu_2}\times\dots\times D_{\mu_m}\times B_{\rho},$$ and
any subsystem of $\Phi$ is equivalent under the action of $W$ to some system of such form.

In the case when $p=2$ subsystems $\Phi_1$ of $\Phi$ are constructed in the following way (see
\cite[Section 3]{Carter}). Let $\lambda$, $\mu$ and $\rho$ be partitions such that
$|\lambda|+|\mu|+|\rho|=n$ with no part of $\mu$ equal to $1$. Let
$\lambda=(\lambda_1,\lambda_2,\dots,\lambda_l)$, $\mu=(\mu_1,\mu_2,\dots,\mu_m)$ and
$\rho=(\rho_1,\rho_2,\dots,\rho_r)$. Let $I_1, I_2,\dots, I_l,$ $J_1, J_2,\dots, J_m,$ $K_1,
K_2,\dots, K_r$ be pairwise non-intersecting subsets of the set $\{1, 2, \dots, n\}$ such that
$|I_a|=\lambda_a$ for $1\leqslant a\leqslant l$, $|J_b|=\mu_b$ for $1\leqslant b\leqslant m$ and
$|K_c|=\rho_c$ for $1\leqslant c\leqslant r$. Define the system $\Phi_1$ to be the union
$$\bigcup_{1\leqslant a\leqslant l}\{e_i - e_j, i\neq j, i,j\in I_\alpha\}\cup\bigcup_{1\leqslant b\leqslant m}
\{\pm e_i \pm e_j, i\neq j, i,j\in J_\beta\}\cup$$$$\cup\bigcup_{1\leqslant c\leqslant r}\{\pm e_i
\pm e_j, \pm e_i, i\neq j, i,j\in K_\gamma\}.$$ Then $\Phi_1$ is a subsystem of $\Phi$ of type
$$A_{\lambda_1-1}\times \dots\times A_{\lambda_l-1}\times D_{\mu_1}\times\dots\times
D_{\mu_m}\times B_{\rho_1}\times\dots\times B_{\rho_r},$$ and any subsystem of $\Phi$ is equivalent
under the action of $W$ to some system of such form.

Subsystems of $\Phi_{\overline H}$ can be obtained in the following way (see \cite[Section
3]{Carter}). Let $\lambda$ and $\mu$ be partitions satisfying the condition $|\lambda|+|\mu|=n$
with no part of $\mu$ equal to $1$. Let $\lambda=(\lambda_1,\lambda_2,\dots,\lambda_l)$ and
$\mu=(\mu_1,\mu_2,\dots,\mu_m)$. Let $I_1,I_2,\dots, I_l, J_1,J_2,\dots, J_m$be pairwise
non-intersecting subsets of the set $\{1,2,\dots,n\}$ such that $|I_a|=\lambda_a$ for all
$1\leqslant a\leqslant l$ and $|J_b|=\mu_b$ for all $1\leqslant b\leqslant m$. Define the system
$\Phi_1$ to be the union
$$\bigcup_{1\leqslant a\leqslant l}\{e_i-e_j: i\neq j, i,j\in I_a\}\cup\bigcup_{1\leqslant b\leqslant m}\{\pm
e_i\pm e_j: i\neq j, i,j\in J_b\}.$$ Then $\Phi_1$ is a subsystem of $\Phi_{\overline H}$ of type
$$A_{\lambda_1-1}\times A_{\lambda_2-1}\times\dots\times A_{\lambda_l-1}\times
D_{\mu_1}\times D_{\mu_2}\times\dots\times D_{\mu_m},$$  and any subsystem of $\Phi$ is equivalent
under the action of $W$ to some system of such form. This implies that subsystems of
$\Phi_{\overline H}$ are also subsystems of $\Phi$. Therefore, reductive subgroups of maximal rank
of $\overline H$ are reductive subgroups of maximal rank of $\overline G$.

For subsystem $\Phi_1$ of $\Phi$ define the set of indices $I_{\Phi_1}$ as follows: $i$ and $-i$,
where $i\in\{1,2,\dots,n\}$, lie in $I_{\Phi_1}$ if and only if $\Phi_1$ contains a root non
orthogonal to $e_i$. Denote by $\overline T_{\Phi_1}$ the subgroup of $\overline T$ consisting of
all diagonal matrices which have $1$'s at places with numbers not in  $I_{\Phi_1}$ along the
diagonal. Denote the group $\langle \overline T_{\Phi_1}, \overline U_r, r\in \Phi_1\rangle$ by
$\overline G_{\Phi_1}$. Let $\Phi_1$ be an indecomposable subsystem and $k=|I_{\Phi_1}|/2$.
Repeating arguments from the previous section, we deduce that if $\Phi_1$ has type $A_{k-1}$, then
$\overline G_{\Phi_1}$ is isomorphic to $GL_k(\overline F_p)$. If $\Phi_1$ has type $D_k$ or $B_k$,
then, as we mentioned above, the root subgroups corresponding to the root from subsystem $\Phi_1$
generate a subgroup isomorphic to $\Omega_{2k}(\overline F_p)$ or $\Omega_{2k+1}(\overline F_p)$
respectively. Let $\overline S$ be a maximal torus of $\Omega_l(\overline F_p)$. If $p$ is odd,
then group $SO_l(\overline F_p)$ is generated by $\overline S$ and $\Omega_l(\overline F_p)$. If
$p=2$, then $\overline S$ is contained in $\Omega_l(\overline F_p)$. Thus if $p$ is odd, then
$\overline G_{\Phi_1}$ is isomorphic to $SO_{2k}(\overline F_p)$, if $\Phi_1$ has type $D_k$, and
isomorphic to $SO_{2k+1}(\overline F_p)$, if $\Phi_1$ has type $B_k$; and if $p=2$, then $\overline
G_{\Phi_1}$ is isomorphic to $\Omega_{2k}(\overline F_p)$, if $\Phi_1$ has type $D_k$, and
isomorphic to $\Omega_{2k+1}(\overline F_p)$, if $\Phi_1$ has type $B_k$.

Let $\overline G_1$ be a $\sigma$-stable reductive subgroup of maximal rank of $\overline G$ with
the root system $\Phi_1$ of the form described above. Denote by $n_{i}$ the number of sets of
cardinality $i$ among the sets $I_a$, by $m_j$ the number of sets of cardinality $j$ among the sets
$J_b$. If $p$ is odd, then $\overline G_1$ is isomorphic to the direct product
$$\prod_i GL_i(\overline F_p)^{n_i}\times \prod_j SO_{2j}(\overline F_p)^{m_j}\times SO_{2\rho+1}(\overline F_p).$$
If $p=2$, then $\overline G_1$ is isomorphic to $$\prod_i GL_i(\overline F_p)^{n_i}\times \prod_j
\Omega_{2j}(\overline F_p)^{m_j}\times\prod_s \Omega_{2s+1}(\overline F_p)^{k_s},$$ where $k_s$ is
the number of sets of cardinality $s$ among the sets $K_c$.

Put $P_A=\prod_i  Sl_{n_i}$, $P_D=\prod_j Sl_{m_j}$. Let $P_B=1$, if $p$ is odd, and $P_B=\prod_s
Sym_{k_s}$ if $p=2$. Then $N_W(W_1)/W_1$ is isomorphic to the direct product $P_A\times P_D\times
P_B$ (see \cite[Section 3]{Carter}).

Recall that to describe the structure of reductive subgroups of maximal rank in finite group it
suffices to describe the structure of groups $(\overline G_1)_{\sigma\circ n}$, where $n$ runs over
some set of matrices. Let us show that matrices $n$ can be chosen among permutation matrices. Let
$w\in N_W(W_1)$. Let $n$  be the permutation matrix corresponding to permutation $w$. Then one of
the matrices $n$ and $-n$ lies in the preimage of $w$. Since actions of $n$ and $-n$ on $\overline
G_1$ coincide, every group of the class corresponding to $W_1w$ is conjugate to $(\overline
G_1)_{\sigma\circ n}$ in $\overline G$.

Let us construct a map from $N_W(W_1)/W_1$ to $N_W(W_1)$. If $\overline w$ is a cycle of length $k$
in a subgroup $Sl_{n_i}$ of $P_A$, then $\overline G_1$ has a direct factor isomorphic to a direct
power $GL_i(\overline F_p)^k$ consisting, after suitable renumbering of rows and columns, of
block-diagonal matrices of the form
$$\operatorname{diag}(X_1, \dots,X_k, E_{n-ki},X_1^{-T},\dots,X_k^{-T},E_{n-ki}),$$ where $X_j\in
GL_i(\overline F_p)$ for $1\leqslant j\leqslant k$. If cycle $\overline w$ is positive, then it
maps to the product of cycles $(1,i+1, \dots, (k-1)i+1)(2,i+2,\dots,(k-1)i+2)\dots(i,2i, \dots,
ki)$. If $\overline w$ is negative, then it maps to the product of cycles $\overline{(1,i+1, \dots,
(k-1)i+1)}\dots\overline{(i,2i, \dots, ki)}$. Let $\overline w$ be a cycle of length $k$ ina
subgroup $Sl_{m_i}$ of $P_D$. Then $\overline G_1$ has a direct factor isomorphic to a direct power
$SO_{2i}(\overline F_p)^k$ consisting, after suitable renumbering of rows and columns, of
block-diagonal matrices of the form $\operatorname{diag}(X_1, X_2,\dots,X_k, E_{2n-2ki})$, $X_j\in
SO_{2i}(\overline F_p)$ for $1\leqslant j\leqslant k$. If $\overline w$ is positive, then it maps
to the same permutation as in the case of a positive cycle of $P_A$. If $\overline w$ is negative,
then it maps to the permutation $(1,i+1,\dots, (k-1)i+1)(2,i+2,\dots,
(k-1)i+2)\dots\overline{(i,2i,\dots, ki)}$. Let $p=2$ and $\overline w$ be a cycle of length $k$ in
a subgroup $Sym_{k_i}$ of $P_B$. Then $\overline G_1$ contains a subgroup isomorphic to
$SO_{2i+1}(\overline F_p)^k$ consisting, after suitable renumbering of rows and columns, of
block-diagonal matrices of the form $\operatorname{diag}(1,X_1, X_2,\dots,X_k, E_{2n-2ki})$, where
$\operatorname{diag}(1,X_j)\in SO_{2i+1}(\overline F_p)$ for all $1\leqslant j\leqslant k$. The
image of $\overline w$ is determined in the same way as in the case of a positive cycle of $P_A$.
It is not difficult to check that the extension of this map to $N_W(W_1)/W_1$ is an isomorphism.
Denote the image of $N_W(W_1)/W_1$ under this isomorphism by $W_2$. Obviously, the intersection of
$W_1$ and $W_2$ is trivial. Thus, $N_W(W_1)=W_1\leftthreetimes W_2$. If we now map each element of
$W_2$ to corresponding permutation matrix of size $(2n+1)\times(2n+1)$, then we will receive an
isomorphic embedding of $N_W(W_1)/W_1$ into the group of permutation matrices. Denote this
embedding by $\phi$. If $W_1$ is a subgroup of $W_{\overline H}$, then group $N_{W_{\overline
H}}(W_1)/W_1$ corresponds to the subgroup $W_2\cap W_{\overline H}$ of $W_2$ consisting of all
permutations of $W_2$,whose decomposition into disjoint cycles contains even number of negative
cycles.

The conjugacy classes of $N_W(W_1)/W_1$ can be determined by collections of partitions. Indeed, a
conjugacy class of $Sym_{k_i}$ is determined by a partition $\rho^{(i)}$ of $k_i$; a conjugacy
class of $Sl_{n_i}$ is determined by a pair of partitions $\xi^{(i)}$ and $\zeta^{(i)}$ with
$|\xi^{(i)}|+|\zeta^{(i)}|=n_i$, and a conjugacy class of $Sl_{m_i}$ is determined by a pair of
partitions $\theta^{(i)}$ and $\upsilon^{(i)}$ with $|\theta^{(i)}|+|\upsilon^{(i)}|=m_i$, where
parts of partitions $\xi^{(i)}$ and $\theta^{(i)}$ give lengths of positive cycles and parts of
$\zeta^{(i)}$ and $\upsilon^{(i)}$ give lengths of negative cycles. We have

\begin{equation}\label{equation2}\sum_i i(|\xi^{(i)}|+|\zeta^{(i)}|)+\sum_j j|\rho^{(i)}|+\sum_s s(|\theta^{(s)}|+|\upsilon^{(s)}|)+\rho=n.
\end{equation}
In this equality $\rho=0$, if the characteristic is $2$, every $|\rho^{(i)}|$ is equal to $0$, if
the characteristic is odd, and for classes of $N_{W_{\overline H}}(W_1)/W_1$ both conditions are
satisfied, that is $\rho$ is equal to $0$ and every $|\rho^{(i)}|$ is equal to $0$.

Given a collection of partitions let $\overline w$ be a representative of the corresponding
conjugacy class. The element $\overline w$ can be presented as a product of disjoint cycles, and
$\overline G_1$ can be presented as a direct product of subgroups such that each cycle acts
nontrivially exactly on one factor. Thus, to describe the structure of reductive subgroups it
suffices, for every cycle, to describe the structure of corresponding subgroup. Group
$N_W(W_1)/W_1$ contains five types of cycles: positive and negative cycles from factors $P_A$,
positive and negative cycles from factors $P_D$ and cycles from factors $P_B$. The first two cases
can be dealt in the same way as in the previous section. The cases of a positive cycle of $P_D$ and
a cycle of $P_B$ can be done in a similar way. Let us consider the case of negative cycle of~$P_D$.

Let $p$ be odd and $\upsilon^{(i)}_j$ be a part of partition $\upsilon^{(i)}$. Then the
decomposition of $\overline w$ into disjoint cycles contains a negative cycle of length
$\upsilon^{(i)}_j$. This implies that $\overline G_1$ contains a subgroup isomorphic to
$SO_{2i}(\overline F_p)^{\upsilon^{(i)}_j}$ consisting of block-diagonal matrices of the form
$$\operatorname{diag}(X_1,X_2 ,\dots, X_{\upsilon^{(i)}_j}, E),$$
where $X_k\in SO_{2i}(\overline F_p)$ for $1\leqslant k\leqslant \upsilon^{(i)}_j$, and $E$ is an
identity matrix of suitable size. We have
\begin{equation*}\sigma\circ\phi(\overline w)(\operatorname{diag}(X_1,X_2 ,\dots,
X_{\upsilon^{(i)}_j}, E))=\operatorname{diag}(X_1,X_2 ,\dots, X_{\upsilon^{(i)}_j},
E).\end{equation*} Rewrite this in the following way \begin{equation*}
\sigma(\operatorname{diag}(X_{\upsilon^{(i)}_j}^\tau,X_1 ,\dots, X_{\upsilon^{(i)}_j-1},
E))=\operatorname{diag}(X_1,X_2 ,\dots, X_{\upsilon^{(i)}_j}, E),
\end{equation*}
where $\tau$ denote the graph automorphism of $SO_{2i}(\overline F_p)$. We obtain the system of
equalities $X_2=X_1^\sigma$, $X_3=X_2^\sigma, \dots, X_1=X_{\upsilon_j^{(i)}}^{\sigma\circ\tau}$.
This system gives equality
$$X_1=X_1^{\sigma^{\upsilon^{(i)}_j}\circ\tau},$$ which is satisfied if and only if the matrix $X_1$ lies
in $SO_{2i}^-(q^{\upsilon^{(i)}_j})$. Thus, in matrices of group $(\overline G_1)_{\sigma\circ
\phi(\overline w)}$, parts $\upsilon^{(i)}_j$ correspond to blocks of the form
$$\operatorname{diag}(X, X^\sigma,\dots, X^{\sigma^{(\upsilon^{(i)}_j-1)}}),\; X\in
SO_{2i}^-(q^{\upsilon^{(i)}_j}).$$ Denote the group of matrices of such from by $H_{ij}^\upsilon$.

In the case $p=2$, by repeating the arguments from the previous paragraph with substitution of
$SO_{2i}(\overline F_p)$ for $\Omega_{2i}(\overline F_p)$, we obtain that, in matrices of
$(\overline G_1)_{\sigma\circ \phi(\overline w)}$, part $\upsilon^{(i)}_j$ corresponds to the block
$$\operatorname{diag}(X, X^\sigma,\dots, X^{\sigma^{(\upsilon^{(i)}_j-1)}}),\; X\in
\Omega_{2i}^-(q^{\upsilon^{(i)}_j}).$$ Denote the group of matrices of such form by
$H_{ij}^\upsilon$.

Groups corresponding to parts of partitions $\xi^{(i)}$ and $\zeta^{(i)}$ can be defined
independently of the characteristic. Define for a part $\xi^{(i)}_j$ of partition $\xi^{(i)}$ the
group $H_{ij}^\xi$ to be the group consisting of matrices of the form $\operatorname{diag}(Y,
Y^{-\top})$, where
$$Y=\operatorname{diag}(X, X^\sigma, X^{\sigma^2},\dots, X^{\sigma^{(\xi^{(i)}_j-1)}}),\;X\in
GL_{i}(q^{\xi^{(i)}_j}).$$ Define for a part $\zeta^{(i)}_j$ of partition $\zeta^{(i)}$ the group
$H_{ij}^\zeta$ to be the group consisting of matrices of the form $\operatorname{diag}(Y,
Y^{-\top})$, where
$$Y=\operatorname{diag}(X, X^\sigma, X^{\sigma^2},\dots, X^{\sigma^{(\zeta^{(i)}_j-1)}}),\;X\in
GU_{i}(q^{\zeta^{(i)}_j}).$$

Let $p$ be odd. Define for a part $\theta^{(i)}_j$ of partition $\theta^{(i)}$ the group
$H_{ij}^\theta$ to be the group consisting of matrices of the form $$\operatorname{diag}(X,
X^\sigma, X^{\sigma^2},\dots, X^{\sigma^{(\theta^{(i)}_j-1)}}), \; X\in
SO^+_{2i}(q^{\theta^{(i)}_j}).$$ If $p=2$, then a part $\theta^{(i)}_j$ of partition $\theta^{(i)}$
corresponds to the group $H_{ij}^\theta$ consisting of matrices of the form
$$\operatorname{diag}(X, X^\sigma, X^{\sigma^2},\dots, X^{\sigma^{(\theta^{(i)}_j-1)}}), \; X\in
\Omega^+_{2i}(q^{\theta^{(i)}_j}).$$

If $p$ is odd, then subsystem of type $B$ is unique and the subgroup corresponding to this
subsystem is isomorphic to $SO_{2\rho+1}$, where $\rho$ satisfies the equality (\ref{equation2}).
Let $p=2$. Define for a part $\rho^{(i)}_j$ of partition $\rho^{(i)}$ the group $H_{ij}^\rho$ to be
the group $$\operatorname{diag}(X, X^\sigma, X^{\sigma^2},\dots,
X^{\sigma^{(\rho^{(i)}_j-1)}}),\;\operatorname{diag}(1,X)\in \Omega_{2i+1}(q^{\rho^{(i)}_j}).$$

Thus, in this notations the following statement, giving the description of reductive subgroups of
$SO_{2n+1}(q)$, holds.

\begin{prop}\label{p:sofin} Let $p$ be odd and $\overline G_1^g$ be $\sigma$-stable. Let $w=\pi(g^\sigma
g^{-1})$ and the coset $W_1w$ lie in the conjugacy class corresponding to a collection of
partitions $\theta^{(i)}$, $\upsilon^{(i)}$, $\xi^{(i)}$, $\zeta^{(i)}$. Then $(\overline
G_1^g)_\sigma$ is conjugate in $\overline G$ to the group $G_1$ consisting of all block-diagonal
matrices with blocks $B$, $A^+_{ij}$, $A^-_{ij}$, $D^+_{ij}$, $D^-_{ij}$, where $B\in
SO_{2\rho+1}(q)$ and $\rho$ satisfies the equality $(\ref{equation2})$, $A^+_{ij}\in H_{ij}^\xi$,
$A^-_{ij}\in H_{ij}^\zeta$, $D^+_{ij}\in H_{ij}^\theta$, $D^-_{ij}\in H_{ij}^\upsilon$. In
particular,
$$\overline G_1^g\simeq
SO_{2\rho+1}(q)\times\prod_{i,j}GL_i(q^{\xi^{(i)}_j})\times\prod_{i,j}GU_i(q^{\zeta^{(i)}_j})\times
\prod_{i,j}SO^+_{2i}(q^{\theta^{(i)}_j})\times\prod_{i,j}SO_{2i}^-(q^{\upsilon^{(i)}_j}).$$
\end{prop}

\begin{prop}\label{p:sofin2} Let $p=2$ and $\overline G_1^g$ be $\sigma$-stable. Let $w=\pi(g^\sigma
g^{-1})$ and the coset $W_1w$ lie in the conjugacy class corresponding to a collection of
partitions $\rho^{(i)},\theta^{(i)}$, $\upsilon^{(i)}$ $\xi^{(i)}$, $\zeta^{(i)}$. Then $(\overline
G_1^g)_\sigma$ is conjugate in $\overline G$ to the group $G_1$ consisting of all block-diagonal
matrices with blocks $B_{ij}$, $A^+_{ij}$, $A^-_{ij}$, $D^+_{ij}$, $D^-_{ij}$, where $B_{ij}\in
H_{ij}^\rho$, $A^+_{ij}\in H_{ij}^\xi$, $A^-_{ij}\in H_{ij}^\zeta$, $D^+_{ij}\in H_{ij}^\theta$,
$D^-_{ij}\in H_{ij}^\upsilon$. In particular, $$\overline G_1^g\simeq
\prod_{i,j}\Omega_{2i+1}(q^{\rho^{(i)}_j})\times\prod_{i,j}GL_i(q^{\xi^{(i)}_j})\times\prod_{i,j}GU_i(q^{\zeta^{(i)}_j})\times
\prod_{i,j}\Omega^+_{2i}(q^{\theta^{(i)}_j})\times\prod_{i,j}\Omega_{2i}^-(q^{\upsilon^{(i)}_j}).$$
\end{prop}

The following proposition is a corollary of Propositions \ref{p:sofin} and \ref{p:sofin2} and gives
a description of reductive subgroups of $SO^\varepsilon_{2i}(q)$ with odd $q$ and
$\Omega^\varepsilon_{2i}(q)$ with $q$ equal to a power of $2$.

\begin{prop}\label{p:so+-fin} Let $\overline G_1^g$ be a $\sigma$-stable subgroup of $\overline H$.
Let $w=\pi(g^\sigma g^{-1})$ and the coset $W_1w$ lie in the conjugacy class corresponding to a
collection of partitions $\theta^{(i)}$, $\upsilon^{(i)}$ $\xi^{(i)}$, $\zeta^{(i)}$. Then
$(\overline G_1^g)_\sigma$ is conjugate in $\overline G$ to the group $G_1$consisting of all
block-diagonal matrices with blocks $A^+_{ij}$, $A^-_{ij}$, $D^+_{ij}$, $D^-_{ij}$, where
$A^+_{ij}\in H_{ij}^\xi$, $A^-_{ij}\in H_{ij}^\zeta$, $D^+_{ij}\in H_{ij}^\theta$, $D^-_{ij}\in
H_{ij}^\upsilon$. Moreover, the number of groups $H_{ij}^\zeta$ with odd $i$ and groups
$H_{ij}^\upsilon$ is even, if $W_1w\subseteq W_{\overline H}$, and odd, if $W_1w\subseteq
w_0W_{\overline H}$. In particular, group $\overline G_1^g$ is isomorphic to
$$\prod_{i,j}GL_i(q^{\xi^{(i)}_j})\times\prod_{i,j}GU_i(q^{\zeta^{(i)}_j})\times
\prod_{i,j}SO^+_{2i}(q^{\theta^{(i)}_j})\times\prod_{i,j}SO_{2i}^-(q^{\upsilon^{(i)}_j}),$$ if $p$
is odd; and isomorphic to
$$\prod_{i,j}GL_i(q^{\xi^{(i)}_j})\times\prod_{i,j}GU_i(q^{\zeta^{(i)}_j})\times
\prod_{i,j}\Omega^+_{2i}(q^{\theta^{(i)}_j})\times\prod_{i,j}\Omega_{2i}^-(q^{\upsilon^{(i)}_j}),$$
in the case $p=2$.
\end{prop}

\begin{prf} Recall that $N_{W}(W_1)/W_1$ can be isomorphically embedded into $N_W(W_1)$. The image
of this embedding we denoted by $W_2$. Under this embedding, positive cycles of $N_W(W_1)/W_1$ map
to products of positive cycles of $N_W(W_1)$. A negative cycle corresponding to a part of partition
$\upsilon^{(i)}$ map to a product of a number of positive cycles and a negative cycle. Finally, a
cycle corresponding to a part of partition $\zeta^{(i)}$ maps to a product of $i$ negative cycles.
To finish the proof we should recall that the image of $N_{W_{\overline H}}(W_1)/W_1$ under the
embedding consists of permutations whose decomposition into disjoint cycles contains an even number
of negative cycles.
\end{prf}

The center of $H_{ij}^\xi$ is isomorphic to a cyclic group of order $q^{\xi^{(i)}_j}-1$, the center
of $H_{ij}^\zeta$ is isomorphic to a cyclic group of order $q^{\zeta^{(i)}_j}+1$. The center of
$H_{ij}^\theta$ is trivial, if $q$ is a power of $2$, and has order $(4,q^{i\theta^{(i)}_j}-1)/2$,
if $q$ is odd. The center of $H_{ij}^\upsilon$ is trivial, if $q$ is a power of $2$, and has order
$(4,q^{i\upsilon^{(i)}_j}+1)/2$, if $q$ is odd. The center of $H_{ij}^\rho$ is trivial. Thus, the
center of $G_1$ is isomorphic to the direct product
\begin{equation}\label{zorth}\prod_{i,j}(q^{\xi^{(i)}_j}-1)\times\prod_{i,j}(q^{\zeta^{(i)}_j}+1)\times
\prod_{i,j}Z(H_{ij}^\theta)\times\prod_{i,j}Z(H_{ij}^\upsilon).\end{equation}

In the next section we will need some information about maximal tori of orthogonal groups with odd
$p$. In this case $\Omega_{n}^\varepsilon(q)=O^{p'}(SO_n^\varepsilon(q))$. For pair of partitions
$\alpha=(\alpha_1,\alpha_2,\dots,\alpha_a)$ and $\beta=(\beta_1,\beta_2,\dots,\beta_b)$ define
elements $t_i$, $1\leqslant i\leqslant a+b$, in the following way: for $1\leqslant i\leqslant a$
put
$$A_i=\operatorname{diag}(E_{\alpha_1+\alpha_2+\dots+\alpha_{i-1}},\lambda_i,\lambda_i^q,\dots,\lambda_i^{q^{\alpha_i-1}},
E_{\alpha_{i+1}+\dots+\alpha_a+|\beta|})$$ and $t_i=\operatorname{diag}(A_i, A_i^{-1}),$ where
$\lambda_i$ is a primitive root of unity of degree $q^{\alpha_i}-1$, for $a< i\leqslant a+b$ put
$$A_i=\operatorname{diag}(E_{|\alpha|+\beta_1+\beta_2+\dots+\beta_{i-a-1}},\lambda_i,\lambda_i^q,\dots,\lambda_i^{q^{\beta_{i-a}-1}},
E_{\beta{i-a+1}+\dots+\beta_b})$$ and $t_i=\operatorname{diag}(A_i, A_i^{-1}),$ where $\lambda_i$
is a primitive root of unity of degree$q^{\beta_{i-a}}+1$.

In the following lemma $\sigma_+$ denotes the map $\sigma$, and $\sigma_-$ denotes the map
$\sigma\circ n_0$.

\begin{lemma}\label{l:toriOmega2n} Let $p$ be odd and $T$ be a maximal torus of $\overline H_{\sigma_\varepsilon}$,
$\varepsilon\in\{+,-\}$ corresponding to a pair of partitions
$\alpha=(\alpha_1,\alpha_2,\dots,\alpha_a)$ and $\beta=(\beta_1,\beta_2,\dots,\beta_b)$, where
parts of $\alpha$ give lengths of positive cycles, and parts of $\beta$ give lengths of negative
cycles. Then group $T\cap O^{p'}(\overline H_{\sigma_\varepsilon})$ is conjugate in $\overline H$
to the group
$$T_1=\{t_1^{k_1}t_2^{k_2}\dots t_s^{k_s}| k_1+k_2+\dots+k_s\text{ is even}\}.$$
\end{lemma}

\begin{prf} See. \cite[Proposition 4.3]{our} and the proof of Theorems 5 and 6 in \cite{our}.
\end{prf}

Since groups $T\cap O^{p'}(\overline H_{\sigma_\varepsilon})$ and $T_1$ are conjugate in $\overline
H$, the center $Z(\overline H_{\sigma_\varepsilon})$ is a subgroup of $T_1$. Therefore, the images
of $T\cap O^{p'}(\overline H_{\sigma_\varepsilon})$ and $T_1$ in $\overline H/Z(\overline
H_{\sigma_\varepsilon})$ are conjugate. The generator of the center $z$ is equal to
$t_1^{|t_1|/2}t_2^{|t_2|/2}\dots t_s^{|t_s|/2}=\operatorname{diag}(1,-E)$.

\sect{Spectra of orthogonal groups}

First, we will obtain a description of spectra of simple orthogonal groups over fields of
characteristic $2$.

\begin{theorem}\label{t:omega2} Let $G=\Omega_{2n+1}(q)$, where $n\geqslant 2$ and $q$ is a power of $2$.
Assume that for every natural $k\geqslant 2$ such that $n_0=2^{k-2}+1<n$ and for ever pair of
partitions $\alpha=(\alpha_1, \alpha_2,\dots, \alpha_a)$ and $\beta=(\beta_1, \beta_2,\dots,
\beta_b)$ such that $n-n_0=|\alpha|+|\beta|$ the set $\nu(G)$ contains a number
$$2^k[q^{\alpha_1}-1, q^{\alpha_2}-1, \dots, q^{\alpha_a}-1, q^{\beta_1}+1, q^{\beta_2}+1,\dots, q^{\beta_b}+1];$$
for every pair of partitions $\gamma=(\gamma_1, \gamma_2,\dots, \gamma_c)$ and $\delta=(\delta_1,
\delta_2,\dots, \delta_d)$ such that $n-1=|\gamma|+|\delta|$ the set $\nu(G)$ contains a number
$$2[q^{\gamma_1}-1, q^{\gamma_2}-1, \dots, q^{\gamma_c}-1, q^{\delta_1}+1, q^{\delta_2}+1,\dots, q^{\delta_d}+1];$$
and does not contain any other number. Then $\mu_m(G)\subseteq \nu(G) \subseteq \omega(G)$.
\end{theorem}

\begin{prf} We shall prove the inclusion $\nu(G)\subseteq \omega(G)$ first. To do this we for every possible $k$
and partitions $\alpha$ and $\beta$, and partitions $\gamma$ and $\delta$, shall construct a
reductive subgroup $H$ of $G$ such that $\eta(H)$ is equal to corresponding element of $\nu(G)$.

Consider the first case. Put $I=\{1, 2,\dots, n_0\}$. Let $\Phi_1=\{\pm e_i \pm e_j, \pm e_i, i\neq
j, i,j\in I\}$. Then $\Phi_1$ is a subsystem of system $\Phi$ of type $B_{n_0}$. Let $\overline
G_1=\langle \overline T, \overline U_r, r\in\Phi_1\rangle$. Let $W_1$ be the Weyl group of
$\overline G_1$. The group $N_W(W_1)/W_1$ is isomorphic to $Sl_{n-n_0}$. Consider the conjugacy
class of $N_W(W_1)/W_1$ corresponding to the pair of partitions $\alpha$ and $\beta$, where parts
$\alpha_i$ give lengths of positive cycles and parts $\beta_j$ give lengths of negative cycles.
Group $G$ contains a reductive subgroup of maximal rank $(\overline G_1^g)_\sigma$, whose structure
is determined by this class. Let $a$ be the number of elements in the partition $\alpha$ and $b$ be
the number of element in the partition $\beta$. By Proposition \ref{p:sofin2} we have
$$(\overline G_1^g)_\sigma\simeq \Omega_{2n_0+1}(q)\times(q^{\alpha_1}-1)\times \dots\times
(q^{\alpha_a}-1)\times(q^{\beta_1}+1)\times\dots\times(q^{\beta_b}+1).$$ By Lemma \ref{l:pperoid}
the maximal order of unipotent element in $\Omega_{2n_0+1}(q)$ is equal to $2^k$. Therefore,
$\eta((\overline G_1^g)_\sigma)=2^k[q^{\alpha_1}-1, q^{\alpha_2}-1, \dots, q^{\alpha_a}-1,
q^{\beta_1}+1, q^{\beta_2}+1,\dots, q^{\beta_b}+1]$.

Let $\gamma$ and $\delta$ be a pair of partitions such that $|\gamma|+|\delta|=n-1$. Put
$\Phi_1=\{\pm e_1\}$. Then $\Phi_1$ is a subsystem of system $\Phi$ of type $B_1$. The group
$N_W(W_1)/W_1$ is isomorphic to $Sl_{n-1}$. Consider the conjugacy class of $N_W(W_1)/W_1$
corresponding to the pair of partitions $\gamma$ and $\delta$, where parts of $\gamma_i$ give
lengths of positive cycles and parts of $\delta_j$ give lengths of negative cycles. Group $G$
contains a reductive subgroup of maximal rank $(\overline G_1^g)_\sigma$, whose structure is
determined by this class. Let $c$ be the number of elements in the partition $\gamma$ and $d$ be
the number of element in the partition $\delta$. By Proposition \ref{p:sofin2} we have
$$(\overline G_1^g)_\sigma\simeq \Omega_{3}(q)\times(q^{\gamma_1}-1)\times\dots\times
(q^{\gamma_c}-1)\times(q^{\delta_1}+1)\times\dots\times(q^{\delta_d}+1).$$ By Lemma \ref{l:pperoid}
the maximal order of unipotent element in $\Omega_{3}(q)$ is equal to $2$. Therefore,
$\eta((\overline G_1^g)_\sigma)=2[q^{\gamma_1}-1, q^{\gamma_2}-1, \dots, q^{\gamma_c}-1,
q^{\delta_1}+1, q^{\delta_2}+1,\dots, q^{\delta_d}+1]$.

It remains to prove the inclusion $\mu_m(G)\subseteq \nu(G)$. To do this we shall show that for
every proper reductive subgroup of maximal rank $H$, which is not a maximal torus, there exists an
element of $\nu(G)$ divided by $\eta(H)$. Let the structure of $H$ be determined by a collection of
partitions $\xi^{(i)}$ and $\zeta^{(i)}$, where $1\leqslant i\leqslant a$, $\theta^{(i)}$ and
$\upsilon^{(i)}$, where $2\leqslant i\leqslant b$, $\rho^{(i)}$, where $1\leqslant i\leqslant c$.
Note that $H$ is a maximal torus if and only if $a=1$, $b=0$ and $c=0$. By formula (\ref{zorth})
the center of $H$ is isomorphic to the direct product
$$\prod_{i,j}(q^{\xi^{(i)}_j}-1)\times\prod_{i,j}(q^{\zeta^{(i)}_j}+1).$$ We have
$$\eta(H)=2^k[\underset{i,j}{\operatorname{lcm}}\{q^{\xi^{(i)}_j}-1\},
\underset{i,j}{\operatorname{lcm}} \{q^{\zeta^{(i)}_j}+1\}],$$ where $k$ is such that $2^k$ is the
least power of $2$ which is greater than $\operatorname{max}\{a-1, 2b-3, 2c-1\}$, i.\,e. the
greatest power of $2$ lying in $\omega(H)$. Put $n_0=2^{k-2}+1$, if $k\geqslant 2$, and $n_0=1$, if
$k=1$. Then $n_0\leqslant\max\{(a+1)/2, b-1, c\}$.

By formula (\ref{equation2}) we have
$$|\xi^{(1)}|+|\zeta^{(1)}|+2(|\xi^{(2)}|+|\zeta^{(2)}|)+\dots+a(|\xi^{(a)}|+|\zeta^{(a)}|)+$$
$$+2(|\theta^{(2)}|+|\upsilon^{(2)}|)+\dots+b(|\theta^{(b)}|+|\upsilon^{(b)}|)+|\rho^{(1)}|+2|\rho^{(2)}|+\dots+c|\rho^{(c)}|=n.$$
Put
$$x=|\xi^{(1)}|+|\zeta^{(1)}|+|\xi^{(2)}|+|\zeta^{(2)}|+\dots+|\xi^{(a)}|+|\zeta^{(a)}|+n_0.$$
Let us prove that $x\leqslant n$. By substituting the expressions for $x$ and $n$, we obtain the
inequality
$$|\xi^{(1)}|+|\zeta^{(1)}|+|\xi^{(2)}|+|\zeta^{(2)}|+\dots+|\xi^{(a)}|+|\zeta^{(a)}|+n_0\leqslant$$
$$\leqslant|\xi^{(1)}|+|\zeta^{(1)}|+2(|\xi^{(2)}|+|\zeta^{(2)}|)+\dots+a(|\xi^{(a)}|+|\zeta^{(a)}|)+$$
$$+2(|\theta^{(2)}|+|\upsilon^{(2)}|)+\dots+b(|\theta^{(b)}|+|\upsilon^{(b)}|)+|\rho^{(1)}|+
2|\rho^{(2)}|+\dots+c|\rho^{(c)}|.$$ Assume that $\operatorname{max}\{(a+1)/2, b-1, c\}=(a+1)/2$.
If $a>2$, then
$$n_0+|\xi^{(a)}|+|\zeta^{(a)}|\leqslant\frac{a+1}{2}+|\xi^{(a)}|+|\zeta^{(a)}|\leqslant
a-1+|\xi^{(a)}|+|\zeta^{(a)}| \leqslant a(|\xi^{(a)}|+|\zeta^{(a)}|).$$ If $a=2$, then $n_0=1$. In
this case we have
$$n_0+|\xi^{(a)}|+|\zeta^{(a)}|= a-1+|\xi^{(a)}|+|\zeta^{(a)}| \leqslant a(|\xi^{(a)}|+|\zeta^{(a)}|).$$
If $\operatorname{max}\{(a+1)/2, b-1, c\}=b-1$, then $n_0\leqslant
b(|\theta^{(b)}|+|\upsilon^{(b)}|)$. If $\operatorname{max}\{(a+1)/2, b-1, c\}=c$, then
$n_0\leqslant c|\rho^{(c)}|$. Thus, in every case the inequality $x\leqslant n$ hlds.

Define partitions $\alpha$ and $\beta$ as follows:
$$\alpha=(\xi^{(1)}_1, \xi^{(1)}_2,\dots, \xi^{(2)}_1, \xi^{(2)}_2,\dots, \xi^{(a)}_1, \dots, n-x),$$
$$\beta=(\zeta^{(1)}_1, \zeta^{(1)}_2,\dots, \zeta^{(2)}_1, \zeta^{(2)}_2,\dots, \zeta^{(a)}_1, \dots).$$

Let $\overline G_1=\langle \overline T, \overline U_r, r\in\Phi_1\rangle$, where $\Phi_1$ is a
subsystem of system $\Phi$ of type $B_{n_0}$. The group $N_W(W_1)/W_1$ is isomorphic to
$Sl_{n-n_0}$. Group $N_W(W_1)/W_1$ contains the conjugacy class corresponding to the pair of
partitions $\alpha$ and $\beta$, where parts of partition $\alpha$ give length of positive cycles
and parts of partition $\beta$ give lengths of negative cycles. Group $G$ contains a reductive
subgroup of maximal rank of the form $(\overline G_1^g)_\sigma$, whose structure is determined by
this class. By formula (\ref{zorth}) the center of this group is isomorphic to the direct product
$$\prod_{i,j}(q^{\xi^{(i)}_j}-1)\times\prod_{i,j}(q^{\zeta^{(i)}_j}+1)\times(q^{n-x}-1).$$ Thus, $$\eta((\overline
G_1^g)_\sigma)=2^k[\underset{i,j}{\operatorname{lcm}}\{q^{\xi^{(i)}_j}-1\},
\underset{i,j}{\operatorname{lcm}} \{q^{\zeta^{(i)}_j}+1\}, q^{n-x}-1].$$ This number is divided by
$\eta(H)$ and lies in $\nu(G)$. The theorem is proved.
\end{prf}

Lemma \ref{l:pperoid}, \cite[Theorem 3]{our} and Theorem \ref{t:omega2} yield a description of
spectra of simple orthogonal groups of odd dimension over fields of characteristic $2$.

\begin{cor}\label{c:specomega2}  Let $q$ be a power of $2$ and $G=\Omega_{2n+1}(q)\simeq Sp_{2n}(q)$,
$n\geqslant 2$.  Then $\omega(G)$ consists of all divisors of the following numbers:
\begin{itemize}

\item[$1)$] $[q^{n_1}+\varepsilon_11,q^{n_2}+\varepsilon_21,\dots,q^{n_s}+\varepsilon_s1]$ for all $s\geqslant
1$, $\varepsilon_i\in\{+,-\}$, $1\leqslant i\leqslant s$, and $n_1,n_2,\dots,n_s>0$ such that
$n_1+n_2+\dots+n_s=n$;

\item[$2)$] $2[q^{n_1}+\varepsilon_11,q^{n_2}+\varepsilon_21,\dots,q^{n_s}+\varepsilon_s1]$ for all $s\geqslant 1$,
$\varepsilon_i\in\{+,-\}$, $1\leqslant i\leqslant s$, and $n_1,n_2,\dots,n_s>0$ such that
$n_1+n_2+\dots+n_s=n-1$;

\item[$3)$] $2^k[q^{n_1}+\varepsilon_11,q^{n_2}+\varepsilon_21,\dots,q^{n_s}+\varepsilon_s1]$ for all
$s\geqslant 1$, $k\geqslant 2$, $\varepsilon_i\in\{+,-\}$, $1\leqslant i\leqslant s$, and
$n_1,n_2,\dots,n_s>0$ such that $2^{k-2}+1+n_1+n_2+\dots+n_s=n$;

\item[$4)$] $2^k$, if $2^{k-2}+1=n$ for some $k\geqslant 2$.
\end{itemize}
\end{cor}

\begin{theorem}\label{t:omega+-2} Let $G=\Omega^\varepsilon_{2n}(q)$, where $n\geqslant 4$, $\varepsilon\in\{+,-\}$
and $q$ is a power of $2$. Let $\nu(G)$ to consist of the following numbers:
\begin{itemize}

\item[$1)$] $2^k[q^{\alpha_1}-1, q^{\alpha_2}-1, \dots, q^{\alpha_a}-1, q^{\beta_1}+1, q^{\beta_2}+1,\dots, q^{\beta_b}+1]$
for all $s\geqslant 1$ and $\alpha_1,\dots,\alpha_a,\beta_1\dots,\beta_b>0$ such that
$2^{k-2}+2+\alpha_1+\dots+\alpha_a+\beta_1+\dots+\beta_b=n$;

\item[$2)$] $2[q^{\alpha_1}-1, q^{\alpha_2}-1, \dots, q^{\alpha_a}-1, q^{\beta_1}+1, q^{\beta_2}+1,\dots, q^{\beta_b}+1]$
for all $s\geqslant 1$ and $\alpha_1,\dots,\alpha_a,\beta_1\dots,\beta_b>0$ such that
$2+\alpha_1+\dots+\alpha_a+\beta_1+\dots+\beta_b=n$;

\item[$3)$] $2[q\pm1, q^{\alpha_1}-1, q^{\alpha_2}-1, \dots, q^{\alpha_a}-1, q^{\beta_1}+1, q^{\beta_2}+1,\dots, q^{\beta_b}+1]$
for all $s\geqslant 1$, $l$ is even, if $\varepsilon=+$, and odd, if $\varepsilon=-$, and
$\alpha_1,\dots,\alpha_a,\beta_1\dots,\beta_b>0$ such that
$2+\alpha_1+\dots+\alpha_a+\beta_1+\dots+\beta_b=n$;

\item[$4)$] $4[q-1, q^{\alpha_1}-1, q^{\alpha_2}-1, \dots, q^{\alpha_a}-1, q^{\beta_1}+1, q^{\beta_2}+1,\dots, q^{\beta_b}+1]$
for all $s\geqslant 1$, $l$ is even, if $\varepsilon=+$, and odd, if $\varepsilon=-$, and
$\alpha_1,\dots,\alpha_a,\beta_1\dots,\beta_b>0$ such that
$3+\alpha_1+\dots+\alpha_a+\beta_1+\dots+\beta_b=n$.

\item[$5)$] $4[q+1, q^{\alpha_1}-1, q^{\alpha_2}-1, \dots, q^{\alpha_a}-1, q^{\beta_1}+1, q^{\beta_2}+1,\dots, q^{\beta_b}+1]$
for all $s\geqslant 1$, odd $l$, if $\varepsilon=+$, and even, if $\varepsilon=-$, and
$\alpha_1,\dots,\alpha_a,\beta_1\dots,\beta_b>0$ such that
$3+\alpha_1+\dots+\alpha_a+\beta_1+\dots+\beta_b=n$.
\end{itemize}
Then $\mu_m(G)\subseteq \nu(G) \subseteq \omega(G)$.
\end{theorem}

\begin{prf} Let $\varepsilon=+$. We shall prove the inclusion $\nu(G)\subseteq \omega(G)$ first.
For every element $m$ of $\nu(G)$ we will find a reductive subgroup $H$ such that $m=\eta(H)$.

Let $m=2^k[q^{\alpha_1}-1, q^{\alpha_2}-1, \dots, q^{\alpha_a}-1, q^{\beta_1}+1,
q^{\beta_2}+1,\dots, q^{\beta_b}+1]$, where $b$ is even and $k\geqslant 2$. Put $n_0=2^{k-2}+2$.
Then subgroup $H$ can be chosen as a subgroup isomorphic to the direct product
$$\Omega^+_{2n_0}(q)\times(q^{\alpha_1}-1)\times (q^{\alpha_2}-1)\times\dots\times
(q^{\alpha_a}-1)\times(q^{\beta_1}+1)\times(q^{\beta_2}+1)\times\dots\times(q^{\beta_b}+1).$$

Let $m=2^k[q^{\alpha_1}-1, q^{\alpha_2}-1, \dots, q^{\alpha_a}-1, q^{\beta_1}+1,
q^{\beta_2}+1,\dots, q^{\beta_b}+1]$, where $b$ is odd and $k\geqslant 2$. Put $n_0=2^{k-2}+2$.
Then subgroup $H$ can be chosen as a subgroup isomorphic to the direct product
$$\Omega^-_{2n_0}(q)\times(q^{\alpha_1}-1)\times (q^{\alpha_2}-1)\times\dots\times
(q^{\alpha_a}-1)\times(q^{\beta_1}+1)\times(q^{\beta_2}+1)\times\dots\times(q^{\beta_b}+1).$$

Let $m=2[q^{\alpha_1}-1, q^{\alpha_2}-1, \dots, q^{\alpha_a}-1, q^{\beta_1}+1, q^{\beta_2}+1,\dots,
q^{\beta_b}+1]$, where $b$ is even. Then subgroup $H$ can be chosen as a subgroup isomorphic to the
direct product
$$\Omega^+_{4}(q)\times(q^{\alpha_1}-1)\times (q^{\alpha_2}-1)\times\dots\times
(q^{\alpha_a}-1)\times(q^{\beta_1}+1)\times(q^{\beta_2}+1)\times\dots\times(q^{\beta_b}+1).$$

Let $m=2[q^{\alpha_1}-1, q^{\alpha_2}-1, \dots, q^{\alpha_a}-1, q^{\beta_1}+1, q^{\beta_2}+1,\dots,
q^{\beta_b}+1]$, where $b$ is odd. Then subgroup $H$ can be chosen as a subgroup isomorphic to the
direct product
$$\Omega^-_{4}(q)\times(q^{\alpha_1}-1)\times (q^{\alpha_2}-1)\times\dots\times
(q^{\alpha_a}-1)\times(q^{\beta_1}+1)\times(q^{\beta_2}+1)\times\dots\times(q^{\beta_b}+1).$$

Let $m=2[q-1,q^{\alpha_1}-1, q^{\alpha_2}-1, \dots, q^{\alpha_a}-1, q^{\beta_1}+1,
q^{\beta_2}+1,\dots, q^{\beta_b}+1]$. Then subgroup $H$ can be chosen as a subgroup isomorphic to
the direct product
$$GL_{2}(q)\times(q^{\alpha_1}-1)\times (q^{\alpha_2}-1)\times\dots\times
(q^{\alpha_a}-1)\times(q^{\beta_1}+1)\times(q^{\beta_2}+1)\times\dots\times(q^{\beta_b}+1).$$

Let $m=2[q+1,q^{\alpha_1}-1, q^{\alpha_2}-1, \dots, q^{\alpha_a}-1, q^{\beta_1}+1,
q^{\beta_2}+1,\dots, q^{\beta_b}+1]$. Then subgroup $H$ can be chosen as a subgroup isomorphic to
the direct product
$$GU_{2}(q)\times(q^{\alpha_1}-1)\times (q^{\alpha_2}-1)\times\dots\times
(q^{\alpha_a}-1)\times(q^{\beta_1}+1)\times(q^{\beta_2}+1)\times\dots\times(q^{\beta_b}+1).$$

Let $m=4[q-1,q^{\alpha_1}-1, q^{\alpha_2}-1, \dots, q^{\alpha_a}-1, q^{\beta_1}+1,
q^{\beta_2}+1,\dots, q^{\beta_b}+1]$. Then subgroup $H$ can be chosen as a subgroup isomorphic to
the direct product
$$GL_{3}(q)\times(q^{\alpha_1}-1)\times (q^{\alpha_2}-1)\times\dots\times
(q^{\alpha_a}-1)\times(q^{\beta_1}+1)\times(q^{\beta_2}+1)\times\dots\times(q^{\beta_b}+1).$$

Let $m=4[q+1,q^{\alpha_1}-1, q^{\alpha_2}-1, \dots, q^{\alpha_a}-1, q^{\beta_1}+1,
q^{\beta_2}+1,\dots, q^{\beta_b}+1]$. Then subgroup $H$ can be chosen as a subgroup isomorphic to
the direct product
$$GU_{3}(q)\times(q^{\alpha_1}-1)\times (q^{\alpha_2}-1)\times\dots\times
(q^{\alpha_a}-1)\times(q^{\beta_1}+1)\times(q^{\beta_2}+1)\times\dots\times(q^{\beta_b}+1).$$

Now we will prove the inclusion $\mu_m(G)\subseteq \nu(G)$. Let $H$ be a reductive subgroup of
maximal rank, which is not a maximal torus, whose structure is determined by partitions $\xi^{(i)}$
and $\zeta^{(i)}$, where $1\leqslant i\leqslant a$, $\theta^{(i)}$ and $\upsilon^{(i)}$, where
$2\leqslant i\leqslant b$. Since $H$ is not a maximal torus, either $a>1$ or $b>1$. By formula
(\ref{zorth}) the center of $H$ is isomorphic to the direct product
$$\prod_{i,j}(q^{\xi^{(i)}_j}-1)\times\prod_{i,j}(q^{\zeta^{(i)}_j}+1).$$ We have
$$\eta(H)=2^k[\underset{i,j}{\operatorname{lcm}}\{q^{\xi^{(i)}_j}-1\},
\underset{i,j}{\operatorname{lcm}} \{q^{\zeta^{(i)}_j}+1\}],$$ where $k$ is such that $2^k$ is the
least power of $2$, which is greater than $\operatorname{max}\{a-1, 2b-3\}$, i.\,e. the greatest
power of $2$ lying in  $\omega(H)$. Put $n_0=2^{k-2}+2$, if $k\geqslant 2$, and $n_0=2$, if $k=1$.
Then $n_0\leqslant\max\{(a+3)/2, b\}$.

We have
$$|\xi^{(1)}|+|\zeta^{(1)}|+2(|\xi^{(2)}|+|\zeta^{(2)}|)+\dots+a(|\xi^{(a)}|+|\zeta^{(a)}|)+$$
$$+2(|\theta^{(2)}|+|\upsilon^{(2)}|)+\dots+b(|\theta^{(b)}|+|\upsilon^{(b)}|)
=n.$$ Put
$$x=|\xi^{(1)}|+|\zeta^{(1)}|+|\xi^{(2)}|+|\zeta^{(2)}|+\dots+|\xi^{(a)}|+|\zeta^{(a)}|+n_0.$$ Let
us find the conditions under which the inequality $x\leqslant n$ holds. Assume that
$\operatorname{max}\{(a+3)/2, b\}=b$. Then $n_0\leqslant b(|\theta^{(b)}|+|\upsilon^{(b)}|)$ and,
therefore, $x\leqslant n$. Assume that $\operatorname{max}\{(a+3)/2, b\}=(a+3)/2$. If $a>4$, then
$$n_0+|\xi^{(a)}|+|\zeta^{(a)}|\leqslant\frac{a+3}{2}+|\xi^{(a)}|+|\zeta^{(a)}|\leqslant a-1+|\xi^{(a)}|+|\zeta^{(a)}|
\leqslant a(|\xi^{(a)}|+|\zeta^{(a)}|).$$ If $a=4$, then $n_0=3$. In this case we have
$$n_0+|\xi^{(a)}|+|\zeta^{(a)}|= a-1+|\xi^{(a)}|+|\zeta^{(a)}| \leqslant a(|\xi^{(a)}|+|\zeta^{(a)}|).$$
If $a=2$ or $a=3$, then $x$ can be greater than $n$.

Let $a=3$. Then $n_0=3$ and $b\leqslant 3$. Assume that $x>n$. We have
$$|\xi^{(1)}|+|\zeta^{(1)}|+|\xi^{(2)}|+|\zeta^{(2)}|+|\xi^{(3)}|+|\zeta^{(3)}|+3>
|\xi^{(1)}|+|\zeta^{(1)}|+2(|\xi^{(2)}|+|\zeta^{(2)}|)+3(|\xi^{(3)}|+|\zeta^{(3)}|)+$$
$$+2(|\theta^{(2)}|+|\upsilon^{(2)}|)+3(|\theta^{(3)}|+|\upsilon^{(3)}|).$$
Therefore,
$$3>|\xi^{(2)}|+|\zeta^{(2)}|+2(|\xi^{(3)}|+|\zeta^{(3)}|)+
2(|\theta^{(2)}|+|\upsilon^{(2)}|)+3(|\theta^{(3)}|+|\upsilon^{(3)}|).$$ This implies that
$|\xi^{(3)}|+|\zeta^{(3)}|=1$, $|\xi^{(2)}|+|\zeta^{(2)}|=0$ and $b=0$. The group $H$ corresponding
to such collection of partitions is isomorphic to
$$GL_{3}(q)\times\prod_{j}(q^{\xi^{(1)}_j}-1)\times\prod_{j}(q^{\zeta^{(1)}_j}+1),$$ if
$|\xi^{(3)}|=1$, and isomorphic to
$$GU_{3}(q)\times\prod_{j}(q^{\xi^{(1)}_j}-1)\times\prod_{j}(q^{\zeta^{(1)}_j}+1),$$ if
$|\zeta^{(3)}|=1$. In both cases $\eta(H)$ lies in $\nu(G)$.

Let $a=2$. Then $n_0=2$ and $b\leqslant 2$. Then the inequality $x>n$ implies that
$|\xi^{(2)}|+|\zeta^{(2)}|=1$ and $b=0$. The group $H$ corresponding to such collection of
partitions is isomorphic to
$$GL_{2}(q)\times\prod_{j}(q^{\xi^{(1)}_j}-1)\times\prod_{j}(q^{\zeta^{(1)}_j}+1),$$ if
$|\xi^{(2)}|=1$, and isomorphic to
$$GU_{2}(q)\times\prod_{j}(q^{\xi^{(1)}_j}-1)\times\prod_{j}(q^{\zeta^{(1)}_j}+1),$$ if
$|\zeta^{(2)}|=1$. In both cases $\eta(H)$ lies in $\nu(G)$.

Thus, it remains to consider the cases when $x\leqslant n$. If the sum of lengths of partitions
$\zeta^{(i)}$ over all $i$ is even, then $G$ contains a subgroup $H_1$ isomorphic to the direct
product

$$\Omega^+_{2n_0}(q)\times\prod_{i,j}(q^{\xi^{(i)}_j}-1)\times\prod_{i,j}(q^{\zeta^{(i)}_j}+1)\times(q^{n-x}-1),$$
if the sum is odd, then $G$ contains a subgroup $H_1$ isomorphic to

$$\Omega^-_{2n_0}(q)\times\prod_{i,j}(q^{\xi^{(i)}_j}-1)\times\prod_{i,j}(q^{\zeta^{(i)}_j}+1)\times(q^{n-x}-1).$$
In both cases we have $$\eta(H_1)=2^k[\underset{i,j}{\operatorname{lcm}}\{q^{\xi^{(i)}_j}-1\},
\underset{i,j}{\operatorname{lcm}} \{q^{\zeta^{(i)}_j}+1\}, q^{n-x}-1].$$ This number is divided by
$\eta(H)$ and lies in $\nu(G)$. By word for word repeating of this proof with interchanging the
word "even"\ and "odd"{}, we obtain the proof for the case $\varepsilon=-$. The theorem is proved.

\end{prf}

As a corollary of Lemma \ref{l:pperoid}, \cite[Theorem 7]{our} and Theorem \ref{t:omega+-2} we
obtain a description of spectra of simple orthogonal groups of eve dimension over field of
characteristic $2$.

\begin{cor}\label{c:omega+-2} Let $G=\Omega^\varepsilon_{2n}(q)$, where $n\geqslant 4$, $\varepsilon\in\{+,-\}$
and $q$ is a power of $2$. Then $\omega(G)$ consists of all divisors of the following numbers:
\begin{itemize}

\item[$1)$] $[q^{n_1}+1, q^{n_2}+1, \dots, q^{n_l}+1, q^{n_{l+1}}-1, q^{n_{l+2}}-1,\dots, q^{n_s}-1]$
for all $s\geqslant 1$, $l$ is even, if $\varepsilon=+$, and odd, if $\varepsilon=-$, and
$n_1,n_2,\dots,n_s>0$ such that $n_1+n_2+\dots+n_s=n$;

\item[$2)$] $2^k[q^{n_1}+1, q^{n_2}+1, \dots, q^{n_l}+1, q^{n_{l+1}}-1, q^{n_{l+2}}-1,\dots, q^{n_s}-1]$
for all $s\geqslant 1$ and $n_1,n_2,\dots,n_s>0$ such that $2^{k-2}+2+n_1+n_2+\dots+n_s=n$;

\item[$3)$] $2[q^{n_1}+1, q^{n_2}+1, \dots, q^{n_l}+1, q^{n_{l+1}}-1, q^{n_{l+2}}-1,\dots, q^{n_s}-1]$
for all $s\geqslant 1$ and $n_1,n_2,\dots,n_s>0$ such that $2+n_1+n_2+\dots+n_s=n$;

\item[$4)$] $2[q\pm1, q^{n_1}+1, q^{n_2}+1, \dots, q^{n_l}+1, q^{n_{l+1}}-1, q^{n_{l+2}}-1,\dots, q^{n_s}-1]$
for all $s\geqslant 1$, $l$ is even, if $\varepsilon=+$, and odd, if $\varepsilon=-$, and
$n_1,n_2,\dots,n_s>0$ such that $2+n_1+n_2+\dots+n_s=n$;

\item[$5)$] $4[q-1, q^{n_1}+1, q^{n_2}+1, \dots, q^{n_s}+1]$
for all $s\geqslant 1$ which is even, if $\varepsilon=+$, and odd, if $\varepsilon=-$, and
$n_1,n_2,\dots,n_s>0$ such that $3+n_1+n_2+\dots+n_s=n$;

\item[$6)$] $4[q+1, q^{n_1}+1, q^{n_2}+1, \dots, q^{n_l}+1, q^{n_{l+1}}-1, q^{n_{l+2}}-1,\dots, q^{n_s}-1]$
for all $s\geqslant 1$, $l$ is odd, if $\varepsilon=+$, and even, if $\varepsilon=-$, and
$n_1,n_2,\dots,n_s>0$ such that $3+n_1+n_2+\dots+n_s=n$;

\item[$7)$] $2^k$, if $n=2^{k-2}+2$ for some $k>2$.
\end{itemize}
\end{cor}

Let us proceed to a description of spectra of orthogonal groups over fields of odd characteristic.

\begin{theorem}\label{t:so} Let $G=SO_{2n+1}(q)$, where $n\geqslant 2$ and $q$ is a power of an odd prime $p$.
Assume that for every $k\geqslant 1$ such that $n_0=(p^{k-1}+1)/2<n$ and for every pair of
partitions $\alpha=(\alpha_1, \alpha_2,\dots, \alpha_a)$ and $\beta=(\beta_1, \beta_2,\dots,
\beta_b)$ such that $n-n_0=|\alpha|+|\beta|$ the set $\nu(G)$ contains a number
$$p^k[q^{\alpha_1}-1, q^{\alpha_2}-1, \dots, q^{\alpha_a}-1, q^{\beta_1}+1, q^{\beta_2}+1,\dots, q^{\beta_b}+1],$$
and does not contain any other number. Then $\mu_m(G)\subseteq \nu(G) \subseteq \omega(G)$.
\end{theorem}

\begin{prf} Let us prove the inclusion $\nu(G)\subseteq \omega(G)$. Let $H$ be a reductive subgroup
of $G$ isomorphic to
$$SO_{2n_0+1}(q)\times(q^{\alpha_1}-1)\times (q^{\alpha_2}-1)\times\dots\times
(q^{\alpha_a}-1)\times(q^{\beta_1}+1)\times(q^{\beta_2}+1)\times\dots\times(q^{\beta_b}+1).$$ Then
$\eta(H)=p^k[q^{\alpha_1}-1, q^{\alpha_2}-1, \dots, q^{\alpha_a}-1, q^{\beta_1}+1,
q^{\beta_2}+1,\dots, q^{\beta_b}+1]$.

Let us prove that $\mu_m(G)\subseteq \nu(G)$. Let $H$ be a reductive subgroup of maximal rank of
$G$, which is not a maximal torus. Suppose that the structure of $H$ is determined by partitions
$\xi^{(i)}$ and $\zeta^{(i)}$, where $1\leqslant i\leqslant a$, $\theta^{(i)}$ and
$\upsilon^{(i)}$, where $2\leqslant i\leqslant b$. Since $H$ is not a maximal torus, either $a>1$,
or $b>1$, or $\rho>0$, where $\rho$ satisfies the equality (\ref{equation2}). By formula
(\ref{zorth}) the center of $H$ is isomorphic to the direct product
$$\prod_{i,j}(q^{\xi^{(i)}_j}-1)\times\prod_{i,j}(q^{\zeta^{(i)}_j}+1)\times
\prod_{i,j}Z(H_{ij}^\theta)\times\prod_{i,j}Z(H_{ij}^\upsilon).$$ Put
$e=\exp(\prod_{i,j}Z(H_{ij}^\theta)\times\prod_{i,j}Z(H_{ij}^\upsilon)).$ Note that $e$ can take
only two values: $1$ or $2$. We have
$$\eta(H)=p^k[\underset{i,j}{\operatorname{lcm}}\{q^{\xi^{(i)}_j}-1\},
\underset{i,j}{\operatorname{lcm}} \{q^{\zeta^{(i)}_j}+1\},e],$$ where $k$ is such that $p^k$ is
the least power of $p$ which is greater than $\operatorname{max}\{a-1, 2b-3, 2\rho-1\}$, i.\,e. the
greatest power of $p$ lying in $\omega(H)$. Put $n_0=(p^{k-1}+1)/2$. Then $n_0\leqslant\max\{a/2,
b-1, \rho\}$.

We have
$$|\xi^{(1)}|+|\zeta^{(1)}|+2(|\xi^{(2)}|+|\zeta^{(2)}|)+\dots+a(|\xi^{(a)}|+|\zeta^{(a)}|)+$$
$$+2(|\theta^{(2)}|+|\upsilon^{(2)}|)+\dots+b(|\theta^{(b)}|+|\upsilon^{(b)}|)+\rho=n.$$
Put
$$x=|\xi^{(1)}|+|\zeta^{(1)}|+|\xi^{(2)}|+|\zeta^{(2)}|+\dots+|\xi^{(a)}|+|\zeta^{(a)}|+n_0.$$
Let us show that $x\leqslant n$. Assume that $\operatorname{max}\{a/2, b-1, \rho\}=a/2$. Then
$$n_0+|\xi^{(a)}|+|\zeta^{(a)}|\leqslant\frac{a}{2}+|\xi^{(a)}|+|\zeta^{(a)}|\leqslant a-1+|\xi^{(a)}|+|\zeta^{(a)}|
\leqslant a(|\xi^{(a)}|+|\zeta^{(a)}|).$$ If $\operatorname{max}\{a/2, b-1, \rho\}=b-1$, then
$n_0\leqslant b(|\theta^{(b)}|+|\upsilon^{(b)}|)$. If $\operatorname{max}\{a/2, b-1, \rho\}=\rho$,
then $n_0\leqslant \rho$. Thus, in all cases we obtain that $x\leqslant n$. Therefore, group $G$
contains a reductive subgroup $H_1$ isomorphic to
$$SO_{2n_0+1}\times\prod_{i,j}(q^{\xi^{(i)}_j}-1)\times\prod_{i,j}(q^{\zeta^{(i)}_j}+1)\times(q^{n-x}-1).$$
We have
$$\eta(H_1)=p^k[\underset{i,j}{\operatorname{lcm}}\{q^{\xi^{(i)}_j}-1\},
\underset{i,j}{\operatorname{lcm}} \{q^{\zeta^{(i)}_j}+1\}, q^{n-x}-1].$$ This number is divided by
$\eta(H)$ if ay least one of the following conditions holds: either $a>0$, or $n>x$, or $e=1$.
Assume that $a=0$, $n=x$ and $e=2$. Then $n=x=n_0$ and, therefore, $H=G$, but then $e=1$; a
contradiction. The theorem is proved.
\end{prf}

As a corollary of Lemma \ref{l:pperoid}, \cite[Proposition 4.1]{our} and Theorem \ref{t:so} we
obtain the following statement.

\begin{cor}\label{t:specSO}  Let $G=SO_{2n+1}(q)$, where $n\geqslant 2$, $q$ is a power of an odd prime $p$.
Then $\omega(G)$ consists of all divisors of the following numbers:
\begin{itemize}
\item[$1)$] $[q^{n_1}+\varepsilon_11,q^{n_2}+\varepsilon_21,\dots,q^{n_s}+\varepsilon_s1]$ for all $s\geqslant 1$,
$\varepsilon_i\in\{+,-\}$, $1\leqslant i\leqslant s$, and $n_1,n_2,\dots,n_s>0$ such that
$n_1+n_2+\dots+n_s=n$;

\item[$2)$] $p^k[q^{n_1}+\varepsilon_11,q^{n_2}+\varepsilon_21,\dots,q^{n_s}+\varepsilon_s1]$ for
all $s\geqslant 1$, $\varepsilon_i\in\{+,-\}$, $1\leqslant i\leqslant s$, and
$k,n_1,n_2,\dots,n_s>0$ such that $p^{k-1}+1+2n_1+2n_2+\dots+2n_s=2n$;

\item[$3)$] $p^k$, if $p^{k-1}+1=2n$ for some $k>0$.
\end{itemize}
\end{cor}

We will need the following lemma to obtain a description of spectra of groups $\Omega_{2n+1}(q)$.

\begin{lemma}\label{l:intersectionOmega} Let $H$ be a proper reductive subgroup of maximal rank of
$SO_{2n+1}(q)$, where $n>2$ and $q$ is a power of an odd prime $p$, isomorphic to
$$SO_{2l+1}(q)\times(q^{\alpha_1}-1)\times(q^{\alpha_2}-1)\times\dots\times(q^{\alpha_a}-1)\times(q^{\beta_1}+1)\times
(q^{\beta_2}+1)\times\dots\times(q^{\beta_b}+1),$$ где
$l+\alpha_1+\dots+\alpha_a+\beta_1+\dots+\beta_b=n$. Let $p^k$ be the maximal power of$p$ contained
in the spectrum of $SO_{2l+1}(q)$. Then $\eta(H\cap\Omega_{2n+1}(q))$ is equal to
$$p^k[q^{\alpha_1}-1, q^{\alpha_2}-1, \dots, q^{\alpha_a}-1, q^{\beta_1}+1, q^{\beta_2}+1,\dots, q^{\beta_b}+1],$$
if $a+b\geqslant 2$, equal to $p^k(q^{\alpha_1}-1)/2$, if $a=1$ and $b=0$, and equal to
$p^k(q^{\beta_1}+1)/2$, if $a=0$ and $b=1$.
\end{lemma}

\begin{prf} Let $G_1$ be the subgroup of $SO_{2n+1}(q)$ consisting of all matrices of the form
$\operatorname{diag}(A,E)$, where $A\in SO_{2(n-l)+1}(q)$. The center of $H$ is isomorphic to the
direct product
$(q^{\alpha_1}-1)\times\dots\times(q^{\alpha_a}-1)\times(q^{\beta_1}+1)\times\dots\times(q^{\beta_b}+1)$
and conjugate in $SO_{2n+1}(\overline F_p)$ to some maximal torus $S$ of $G_1$ corresponding to the
pair of partitions $\alpha=(\alpha_1,\alpha_2,\dots,\alpha_a)$ and
$\beta=(\beta_1,\beta_2,\dots,\beta_b)$, where parts of partition $\alpha$ give lengths of positive
cycles and parts of partition $\beta$ give lengths of negative cycles. By \cite[Lemma 4.2]{our} the
center of $H\cap\Omega_{2n+1}(q)$ is isomorphic to the intersection of $S$ and the commutator
subgroup $O^{p'}(G_1)$ of $G_1$. Finally, \cite[Theorem 5]{our} yeilds that $\exp(S\cap
O^{p'}(G_1))=\exp(S)$, if $a+b>1$, and $\exp(S\cap O^{p'}(G_1))=\exp(S)/2$, if $a+b=1$. The lemma
is proved.
\end{prf}

Since $\Omega_{5}(q)\simeq PSp_4(q)$, the case $n=2$ in the following theorem is not considered.

\begin{theorem}\label{t:omega} Let $G=\Omega_{2n+1}(q)$, where $n\geqslant 3$ and $q$ is a power of an odd prime
$p$. Assume that for every natural $k$ such that $2n_0=p^{k-1}+1< 2n$ and for every pair of partitions
$\alpha=(\alpha_1, \alpha_2,\dots, \alpha_a)$ and $\beta=(\beta_1, \beta_2,\dots, \beta_b)$ such that
$n-n_0=|\alpha|+|\beta|$ the set $\nu(G)$ contains the number
$$p^k[q^{\alpha_1}-1, q^{\alpha_2}-1, \dots, q^{\alpha_a}-1, q^{\beta_1}+1, q^{\alpha_2}+1,\dots, q^{\beta_b}+1],$$ if
$a+b\geqslant 2$, contains the number $p^k(q^{\alpha_1}-1)/2$, if $a=1$ and $b=0$, contains the number
$p^k(q^{\beta_1}+1)/2$, if $a=0$ and $b=1$, and does not contain any other number. Then $\mu_m(G)\subseteq \nu(G)
\subseteq \omega(G)$.
\end{theorem}

\begin{prf} Lemma \ref{l:intersectionOmega} implies the inclusion $\nu(G)\subseteq\omega(G)$.
Let us show that $\mu_m(G)\subseteq \nu(G)$. Let $H$ be a reductive subgroup of maximal rank of $G$, which is not a
maximal torus. Consider the subgroup $H_1$ of $SO_{2n+1}(q)$ such that $H=H_1\cap \Omega_{2n+1}(q)$. Suppose that the
structure of $H_1$ is determined by partitions $\xi^{(i)}$, $\zeta^{(i)}$, $1\leqslant i\leqslant a,$ and
$\theta^{(i)}$, $\upsilon^{(i)}$, $2\leqslant i \leqslant b$. Then, as it was mentioned above, its center is isomorphic
to the direct product
$$\prod_{i,j}(q^{\xi^{(i)}_j}-1)\times\prod_{i,j}(q^{\zeta^{(i)}_j}+1)\times
\prod_{i,j}Z(H_{ij}^\theta)\times\prod_{i,j}Z(H_{ij}^\upsilon).$$ Put
$e=\exp(\prod_{i,j}Z(H_{ij}^\theta)\times\prod_{i,j}Z(H_{ij}^\upsilon))$. We have
$$\eta(H_1)=p^k[\underset{i,j}{\operatorname{lcm}}\{q^{\xi^{(i)}_j}-1\},
\underset{i,j}{\operatorname{lcm}} \{q^{\zeta^{(i)}_j}+1\}, e],$$ where $k$ is such that $p^k$ is the least power of
$p$ which is greater than $\operatorname{max}\{a-1, 2b-3, 2\rho-1\}$, i.\,e. the greatest power of $p$ lying in
$\omega(H_1)$. Then $\eta(H)=\eta(H_1)/c$, where $c$ takes values $1$ or $2$ depending on structure of $H_1$. Put
$n_0=(p^{k-1}+1)/2$ and $$x=|\xi^{(1)}|+|\zeta^{(1)}|+|\xi^{(2)}|+|\zeta^{(2)}|+\dots+|\xi^{(a)}|+|\zeta^{(a)}|+n_0.$$
In the proof of Theorem \ref{t:so} we showed that $x\leqslant n.$ Therefore, $SO_{2n+1}(q)$ contains a subgroup $K$
isomorphic to the direct product
$$SO_{2n_0+1}(q)\times\prod_{i,j}(q^{\xi^{(i)}_j}-1)\times\prod_{i,j}(q^{\zeta^{(i)}_j}+1)\times
\underbrace{(q-1)\times(q-1)\times\dots\times(q-1)}_{n-x}.$$ We have
$$\eta(K)=p^k[\underset{i,j}{\operatorname{lcm}}\{q^{\xi^{(i)}_j}-1\},
\underset{i,j}{\operatorname{lcm}} \{q^{\zeta^{(i)}_j}+1\}, b],$$ where $b=q-1$, if $n-x>0$, and $b=1$, if $n-x=0$. If
the center $Z(K)$ is not cyclic, then Lemma \ref{l:intersectionOmega} implies that $\eta(K)=\eta(K\cap
\Omega_{2n+1}(q))$. Thus, in this case $\eta(K\cap \Omega_{2n+1}(q))$ is divided by $\eta(H_1)$ and, therefore, divide
by $\eta(H)$.

The center of $K$ is cyclic only if one of the following conditions holds:

1) $n-x=0$ and the sum of lengths of partitions $\xi^{(i)}$ and $\zeta^{(i)}$ over all $i$ is equal to $1$;

2) $n-x=1$ and $a=0$.

Consider the first case. Recall that
$$n=|\xi^{(1)}|+|\zeta^{(1)}|+2(|\xi^{(2)}|+|\zeta^{(2)}|)+\dots+a(|\xi^{(a)}|+|\zeta^{(a)}|)+$$
$$+2(|\theta^{(2)}|+|\upsilon^{(2)}|)+\dots+b(|\theta^{(b)}|+
|\upsilon^{(b)}|)+\rho.$$

By substituting in equality $n=x$ the expressions for $n$ and $x$, we obtain the following equality
$$n_0=|\xi^{(2)}|+|\zeta^{(2)}|+\dots+(a-1)(|\xi^{(a)}|+|\zeta^{(a)}|)+$$
$$+2(|\theta^{(2)}|+|\upsilon^{(2)}|)+\dots+b(|\theta^{(b)}|+
|\upsilon^{(b)}|)+\rho.$$ Therefore, $n_0\geqslant\max\{a-1,b,\rho\}$. We have $n_0\leqslant\max\{a/2,b-1,\rho\}$.
Assume that $\max\{a/2,b-1,\rho\}=a/2$. Then $a=2$, $n_0=1$ and $b=\rho=0$, i.e. $n=2$ contrary to the hypothesis. The
equality $\max\{a/2,b-1,\rho\}=b-1$ obviously leads to a contradiction. Assume that $\max\{a/2,b-1,\rho\}=\rho$. Then
$n_0=\rho$, $b=0$ and $a=1$. Hence $H_1$ is isomorphic to
$$SO_{2n_0+1}(q)\times\prod_{j}(q^{\xi^{(1)}_j}-1)\times\prod_{j}(q^{\zeta^{(1)}_j}+1).$$
Lemma \ref{l:intersectionOmega} implies that $\eta(H)\in\nu(G)$.

Consider the second case. We have $$n_0=2(|\theta^{(2)}|+|\upsilon^{(2)}|)+\dots+b(|\theta^{(b)}|+
|\upsilon^{(b)}|)+\rho-1.$$ Therefore, $n_0\geqslant\max\{b-1,\rho-1\}$. Thus, $n_0$ is equal to either $b-1$, or
$\rho-1$, or $\rho$. If $n_0$ equals $\rho-1$ or $\rho$, then $b=0$ and $H=G$. Suppose that $n_0=b-1$. Then $\rho=0$,
$n=b$ and $|\theta^{(b)}|+ |\upsilon^{(b)}|=1$. This equalities imply that either $H_1\simeq SO^+_{2n}(q)$, or
$H_1\simeq SO^-_{2n}(q)$. Therefore, $\eta(H)=p^kd$, where $d$ equals $1$ or $2$. Consider reductive subgroups of
maximal rank of $SO_{2n+1}(q)$ isomorphic to $SO_{2n_0+1}(q)\times(q-1)$ and $SO_{2n_0+1}(q)\times(q+1)$. By Lemma
\ref{l:intersectionOmega} the values of $\eta$ on intersection of these subgroups with $\Omega_{2n+1}(q)$ equal to
$p^k(q-1)/2$ and $p^k(q+1)/2$ respectively and lie in $\nu(G)$. One of the numbers $(q-1)/2$ and $(q+1)/2$ is even,
therefore, $\eta(H)$ divides either $p^k(q-1)/2$, or $p^k(q+1)/2$. The theorem is proved.
\end{prf}

As a corollary of Lemma \ref{l:pperoid}, \cite[Theorem 4]{our} and Theorem \ref{t:omega} we obtain a description of
spectra of simple orthogonal groups of odd dimension over fields of odd characteristic.

\begin{cor}\label{c:specOmega}  Let $G=\Omega_{2n+1}(q)$, where $n\geqslant 3$, $q$ is a power of an odd prime
$p$. Then $\omega(G)$ consists of all divisors of the following numbers:
\begin{itemize}
\item[$1)$] $\frac{q^n\pm1}{2}$;

\item[$2)$] $[q^{n_1}+\varepsilon_11,q^{n_2}+\varepsilon_21,\dots,q^{n_s}+\varepsilon_s1]$ for all $s\geqslant 2$,
$\varepsilon_i\in\{+,-\}$, $1\leqslant i\leqslant s$, and $n_1,n_2,\dots,n_s>0$ such that $n_1+n_2+\dots+n_s=n$;

\item[$3)$] $p^k\frac{q^{n_1}\pm1}{2}$ for all $k$ and $n_1$ such that $p^{k-1}+1+2n_1=2n$;

\item[$4)$] $p^k[q^{n_1}+\varepsilon_11,q^{n_2}+\varepsilon_21,\dots,q^{n_s}+\varepsilon_s1]$ for all
$s\geqslant 2$, $\varepsilon_i\in\{+,-\}$, $1\leqslant i\leqslant s$, and $k,n_1,n_2,\dots,n_s>0$ such that
$p^{k-1}+1+2n_1+2n_2+\dots+2n_s=2n$;

\item[$5)$] $p^k$, if $p^{k-1}+1=2n$ for some $k>0$.
\end{itemize}
\end{cor}

\begin{theorem}\label{t:so2n} Let $G=SO^\varepsilon_{2n}(q)$, where $n\geqslant 4$, $\varepsilon\in\{+,-\}$ and $q$ is a power of an odd prime $p$.
Assume that for every natural $k\geqslant 1$ such that $n_0=(p^{k-1}+3)/2<n$ and for every pair of partitions
$\alpha=(\alpha_1, \alpha_2,\dots, \alpha_a)$ and $\beta=(\beta_1, \beta_2,\dots, \beta_b)$ such that
$n-n_0=|\alpha|+|\beta|$ the set $\nu(G)$ contains the number
$$p^k[q^{\alpha_1}-1, q^{\alpha_2}-1, \dots, q^{\alpha_a}-1, q^{\beta_1}+1, q^{\beta_2}+1,\dots, q^{\beta_b}+1],$$
for every pair of partitions $\gamma=(\gamma_1, \gamma_2,\dots, \gamma_c)$ and $\delta=(\delta_1, \delta_2,\dots,
\delta_d)$, where $d$ is even, if $\varepsilon=+$, and odd, if $\varepsilon=-$, such that $n-2=|\gamma|+|\delta|$,
contains the number
$$p[q\pm1 ,q^{\gamma_1}-1, q^{\gamma_2}-1, \dots, q^{\gamma_c}-1, q^{\delta_1}+1, q^{\delta_2}+1,\dots, q^{\delta_d}+1],$$
contains the number $2p^k$, if $2n=p^{k-1}+3$ for some $k>0$, and does not contain any other number. Then
$\mu_m(G)\subseteq \nu(G) \subseteq \omega(G)$.
\end{theorem}

\begin{prf}

We shall prove the inclusion $\nu(G)\subseteq \omega(G)$ first. For every element $m$ of $\nu(G)$ we will find a
reductive subgroup $H$ such that $m=\eta(H)$.

Let $m=p^k[q^{\alpha_1}-1, q^{\alpha_2}-1, \dots, q^{\alpha_a}-1, q^{\beta_1}+1, q^{\beta_2}+1,\dots, q^{\beta_b}+1]$,
where $b$ is even and $|\alpha|+|\beta|=n-n_0$. Then subgroup $H$ can be chosen as a subgroup isomorphic to the direct
product
$$\Omega^+_{2n_0}(q)\times(q^{\alpha_1}-1)\times(q^{\alpha_2}-1)\times\dots\times(q^{\alpha_a}-1)\times(q^{\beta_1}+1)\times
(q^{\beta_2}+1)\times\dots\times(q^{\beta_b}+1).$$

Let $m=p^k[q^{\alpha_1}-1, q^{\alpha_2}-1, \dots, q^{\alpha_a}-1, q^{\beta_1}+1, q^{\beta_2}+1,\dots, q^{\beta_b}+1]$,
where $b$ is odd and $|\alpha|+|\beta|=n-n_0$. Then subgroup $H$ can be chosen as a subgroup isomorphic to the direct
product
$$\Omega^-_{2n_0}(q)\times(q^{\alpha_1}-1)\times(q^{\alpha_2}-1)\times\dots\times(q^{\alpha_a}-1)\times(q^{\beta_1}+1)\times
(q^{\beta_2}+1)\times\dots\times(q^{\beta_b}+1).$$

Let $m=p[q-1 ,q^{\gamma_1}-1, q^{\gamma_2}-1, \dots, q^{\gamma_c}-1, q^{\delta_1}+1, q^{\delta_2}+1,\dots,
q^{\delta_d}+1]$, where $|\gamma|+|\delta|=n-2$. Then subgroup $H$ can be chosen as a subgroup isomorphic to the direct
product
$$GL_2(q)\times(q^{\gamma_1}-1)\times(q^{\gamma_2}-1)\times\dots\times(q^{\gamma_c}-1)\times(q^{\delta_1}+1)\times
(q^{\delta_2}+1)\times\dots\times(q^{\delta_d}+1).$$

Let $m=p[q+1 ,q^{\gamma_1}-1, q^{\gamma_2}-1, \dots, q^{\gamma_c}-1, q^{\delta_1}+1, q^{\delta_2}+1,\dots,
q^{\delta_d}+1]$, где $|\gamma|+|\delta|=n-2$. Then subgroup $H$ can be chosen as a subgroup isomorphic to the direct
product
$$GU_2(q)\times(q^{\gamma_1}-1)\times(q^{\gamma_2}-1)\times\dots\times(q^{\gamma_c}-1)\times(q^{\delta_1}+1)\times
(q^{\delta_2}+1)\times\dots\times(q^{\delta_d}+1).$$

In addition, if $2n=p^{k-1}+3$, then $p^k$ is the maximal power of $p$ contained in $\omega(G)$. The power $p^k$ in not
contained in the spectrum of any proper reductive subgroup. Thus, the number $\eta(G)$, equal to $2p^k$, lies in
$\mu_m(G)$ and, therefore, must lie in $\nu(G)$. If $2n\neq p^{k-1}+3$ for every $k\geqslant 1$, then $\eta(G)$ does
not lie in $\mu(G)$ and is a divisor of $\eta(H)$ for some proper reductive subgroup $H$ of $G$.

Now we shall prove the inclusion $\mu_m(G)\subseteq \nu(G)$. Let $H$ ba a proper reductive subgroup of maximal rank,
which is not a maximal torus, whose structure is determined by partitions $\xi^{(i)}$ and $\zeta^{(i)}$, where
$1\leqslant i\leqslant a$, $\theta^{(i)}$ and $\upsilon^{(i)}$, where $2\leqslant i\leqslant b$. Since $H$ is not a
maximal torus, either $a>1$ or $b>1$. By formula (\ref{zorth}) the center of $H$ is isomorphic to the direct product
$$\prod_{i,j}(q^{\xi^{(i)}_j}-1)\times\prod_{i,j}(q^{\zeta^{(i)}_j}+1)\times
\prod_{i,j}Z(H_{ij}^\theta)\times\prod_{i,j}Z(H_{ij}^\upsilon).$$ Put
$e=\exp(\prod_{i,j}Z(H_{ij}^\theta)\times\prod_{i,j}Z(H_{ij}^\upsilon))$. We have
$$\eta(H)=p^k[\underset{i,j}{\operatorname{lcm}}\{q^{\xi^{(i)}_j}-1\},
\underset{i,j}{\operatorname{lcm}} \{q^{\zeta^{(i)}_j}+1\}, e],$$ where $k$ is such that $p^k$ is the least power of
$p$ which is greater than $\max\{a-1, 2b-3\}$, i.\,e. the greatest power of $p$ lying in $\omega(H)$. Put
$n_0=(p^{k-1}+3)/2$. Then $n_0\leqslant\max\{a/2+1, b\}$.

We have
$$|\xi^{(1)}|+|\zeta^{(1)}|+2(|\xi^{(2)}|+|\zeta^{(2)}|)+\dots+a(|\xi^{(a)}|+|\zeta^{(a)}|)+$$
$$+2(|\theta^{(2)}|+|\upsilon^{(2)}|)+\dots+b(|\theta^{(b)}|+|\upsilon^{(b)}|)
=n.$$

Put
$$x=|\xi^{(1)}|+|\zeta^{(1)}|+|\xi^{(2)}|+|\zeta^{(2)}|+\dots+|\xi^{(a)}|+|\zeta^{(a)}|+n_0.$$ Let
us find the conditions under which the inequality $x\leqslant n$ holds.

Assume that $\max\{a/2+1, b\}=b$. Then $n_0\leqslant b(|\theta^{(b)}|+|\upsilon^{(b)}|)$ and, therefore, $x\leqslant
n$. Assume that $\max\{a/2+1, b\}=a/2+1$. If $a\geqslant4$, then
$$n_0+|\xi^{(a)}|+|\zeta^{(a)}|\leqslant\frac{a}{2}+1+|\xi^{(a)}|+|\zeta^{(a)}|\leqslant a-1+|\xi^{(a)}|+|\zeta^{(a)}|
\leqslant a(|\xi^{(a)}|+|\zeta^{(a)}|).$$ If $a=3$, then $n_0=2$. In this case we have
$$n_0+|\xi^{(a)}|+|\zeta^{(a)}|= a-1+|\xi^{(a)}|+|\zeta^{(a)}| \leqslant a(|\xi^{(a)}|+|\zeta^{(a)}|).$$
If $a=2$, then $x$ can be greater then $n$.

Let $a=2$. Then $n_0=2$ and $b\leqslant 2$. Assume that $x>n$. We have
$$|\xi^{(1)}|+|\zeta^{(1)}|+|\xi^{(2)}|+|\zeta^{(2)}|+2>
|\xi^{(1)}|+|\zeta^{(1)}|+2(|\xi^{(2)}|+|\zeta^{(2)}|)+ 2(|\theta^{(2)}|+|\upsilon^{(2)}|).$$ Therefore,
$$2>|\xi^{(2)}|+|\zeta^{(2)}|+|\theta^{(2)}|+|\upsilon^{(2)}|.$$ This implies that $|\xi^{(2)}|+|\zeta^{(2)}|=1$ and $b=0$.
The group $H$ corresponding to such collection of partitions is isomorphic to
$$GL_2(q)\times\prod_{j}(q^{\xi^{(1)}_j}-1)\times\prod_{j}(q^{\zeta^{(1)}_j}+1),$$ if $|\xi^{(2)}|=1$, and isomorphic
to $$GU_{2}(q)\times\prod_{j}(q^{\xi^{(1)}_j}-1)\times\prod_{j}(q^{\zeta^{(1)}_j}+1),$$ if $|\zeta^{(2)}|=1$. In both
cases $\eta(H)$ lies in $\nu(G)$.

Thus, it remains to consider the cases when $x\leqslant n$. If the sum of lengths of partitions $\zeta^{(i)}$ over all
$i$ is even, then $G$ contains a subgroup $H_1$ isomorphic to the direct product

$$\Omega^+_{2n_0}(q)\times\prod_{i,j}(q^{\xi^{(i)}_j}-1)\times\prod_{i,j}(q^{\zeta^{(i)}_j}+1)\times(q^{n-x}-1),$$
if the sum is odd, then $G$ contains a subgroup $H_1$ isomorphic to

$$\Omega^-_{2n_0}(q)\times\prod_{i,j}(q^{\xi^{(i)}_j}-1)\times\prod_{i,j}(q^{\zeta^{(i)}_j}+1)\times(q^{n-x}-1).$$
In both cases we have $$\eta(H_1)=p^k[\underset{i,j}{\operatorname{lcm}}\{q^{\xi^{(i)}_j}-1\},
\underset{i,j}{\operatorname{lcm}} \{q^{\zeta^{(i)}_j}+1\}, q^{n-x}-1].$$ This number is divided by $\eta(H)$ and lies
in $\nu(G)$. The theorem is proved.
\end{prf}

As a corollary of Lemma \ref{l:pperoid}, \cite[Proposition 4.1]{our} and Theorem \ref{t:so2n} we obtain the following
statement.

\begin{cor}\label{c:so2n} Let $G=SO^\varepsilon_{2n}(q)$, where $n\geqslant 4$, $\varepsilon\in\{+,-\}$ and
$q$ is a power of an odd prime $p$. Then $\omega(G)$ consists of all divisors of the following numbers:
\begin{itemize}

\item[$1)$] $[q^{n_1}+1, q^{n_2}+1, \dots, q^{n_l}+1, q^{n_{l+1}}-1, q^{n_{l+2}}-1,\dots, q^{n_s}-1]$
for all $s\geqslant 1$, $l$ is even, if $\varepsilon=+$, and odd, if $\varepsilon=-$, and $n_1,n_2,\dots,n_s>0$ such
that $n_1+n_2+\dots+n_s=n$;

\item[$2)$] $p^k[q^{n_1}+1, q^{n_2}+1, \dots, q^{n_l}+1, q^{n_{l+1}}-1, q^{n_{l+2}}-1,\dots, q^{n_s}-1]$
for all $s\geqslant 1$ and $n_1,n_2,\dots,n_s>0$ such that $p^{k-1}+3+2n_1+2n_2+\dots+2n_s=2n$;

\item[$3)$] $p[q\pm1, q^{n_1}+1, q^{n_2}+1, \dots, q^{n_l}+1, q^{n_{l+1}}-1, q^{n_{l+2}}-1,\dots, q^{n_s}-1]$
for all $s\geqslant 1$, $l$ is even, if $\varepsilon=+$, and odd, if $\varepsilon=-$, and $n_1,n_2,\dots,n_s>0$ such
that $2+n_1+n_2+\dots+n_s=n$;

\item[$4)$] $2p^k$, if $2n=p^{k-1}+3$ for some $k\geqslant 1$.
\end{itemize}
\end{cor}

In the following lemma $GL^+_2(q)$ denotes the group $GL_2(q)$, and $GL^-_2(q)$ denotes the group $GU_2(q)$.

\begin{lemma}\label{l:mu(SO2n)} $1)$ Let $H$ be a proper reductive subgroup of maximal rank of
$SO^\varepsilon_{2n}(q)$, $\varepsilon\in\{+,-\}$, equal to $H_1\times H_2$, where $H_1$ is isomorphic to
$SO^{\varepsilon_1}_{2k}(q)$, $\varepsilon_1\in\{+,-\}$, $k>1$, and subgroup $H_2$ is isomorphic to the direct product
$$(q^{\alpha_1}-1)\times(q^{\alpha_2}-1)\times\dots\times(q^{\alpha_a}-1)\times(q^{\beta_1}+1)\times
(q^{\beta_2}+1)\times\dots\times(q^{\beta_b}+1),$$ and $k+\alpha_1+\dots+\alpha_a+\beta_2+\dots+\beta_b=n$. Then
$\eta(H\cap\Omega_{2n}^\varepsilon(q))=[2,\eta(H)/2]$, if $a+b=1$ and $Z(\Omega_{2k}^{\varepsilon_1}(q))\neq1$, and
$\eta(H\cap\Omega_{2n}^\varepsilon(q))=\eta(H)$ otherwise.

Let $Z(\Omega_{2n}^\varepsilon(q))\neq 1$ and $\tilde H$ be the image of $H\cap\Omega_{2n}^\varepsilon(q)$ in
$P\Omega_{2n}^\varepsilon(q)$. Then $\eta(\tilde H)=\eta(H)$, if $a+b>1$, and $\eta(\tilde H)=\eta(H)/2$, if $a+b=1$.

$2)$ Let $H$ be a proper reductive subgroup of maximal rank of $SO^\varepsilon_{2n}(q)$, $\varepsilon\in\{+,-\}$, equal
to $H_1\times H_2$, where subgroup $H_1$ is isomorphic to $GL^{\varepsilon_1}_{2}(q)$, $\varepsilon_1\in\{+,-\}$, and
subgroup $H_2$ is isomorphic to the direct product
$$(q^{\alpha_1}-1)\times(q^{\alpha_2}-1)\times\dots\times(q^{\alpha_a}-1)\times(q^{\beta_1}+1)\times
(q^{\beta_2}+1)\times\dots\times(q^{\beta_b}+1),$$ and $2+\alpha_1+\dots+\alpha_a+\beta_2+\dots+\beta_b=n$. Then
$\eta(H\cap\Omega_{2n}^\varepsilon(q))=\eta(H)$, if $a+b>1$, and $\eta(H)=p[q-\varepsilon_11,
(q^{n-2}-\varepsilon1)/2]$, if $a+b=1$.

Let $Z(\Omega_{2n}^\varepsilon(q))\neq 1$ and $\tilde H$ be the image of $H\cap\Omega_{2n}^\varepsilon(q)$ in
$P\Omega_{2n}^\varepsilon(q)$. Then $\eta(\tilde H)=\eta(H\cap\Omega^\varepsilon_{2n}(q))$.
\end{lemma}

\begin{prf} Denote $Z(\Omega_{2n}^\varepsilon(q))$ by $Z$.

$1)$ Let $T_1$ be a maximal torus of $H_1$ isomorphic to the cyclic group of order $(q^k-\varepsilon_11)$. The subgroup
$T$ of $H$, equal to $T_1\times H_2$, is a maximal torus of $SO^{\varepsilon}_{2n}(q)$.

Define a pair of partitions $\alpha$ and $\beta$ as follows: if $\varepsilon=+$, then
$\alpha=(k,\alpha_1,\dots,\alpha_a)$ and $\beta=(\beta_1,\dots,\beta_b)$, if $\varepsilon=-$, then
$\alpha=(\alpha_1,\dots,\alpha_a)$ and $\beta=(k,\beta_1,\dots,\beta_b)$. Then structure of $T$ is determined by the
pair of partitions $\alpha$ and $\beta$, and by Lemma \ref{l:toriOmega2n} the torus $T\cap\Omega_{2n}^\varepsilon(q)$
is conjugate in $\overline H$ to the group $\{t_1^{k_1}t_2^{k_2}\dots t_{a+b+1}^{k_{a+b+1}}|
k_1+k_2+\dots+k_{a+b+1}\text{ is even}\}$, where elements $t_1,\dots,t_{a+b+1}$ are defined before the statement of
Lemma \ref{l:toriOmega2n}. Denote by $t_0$ the element $t_1$, if $\varepsilon=+$, and the element $t_{a+1}$, if
$\varepsilon=-$. Then $T_1$  is conjugate to the group $\langle t_0\rangle$. We have
$$Z(H)=\langle t_0^{|t_0|/2}\rangle\times H_2.$$ Note that the number $|t_0|/2$ is even if and only if
$Z(\Omega_{2k}^\varepsilon(q))\neq 1$.

Let $|t_0|/2$ be even. Then $Z(H\cap \Omega^\varepsilon_{2n}(q))=\langle t_0^{|t_0|/2}\rangle\times (H_2\cap
\Omega^\varepsilon_{2n}(q))$. Now, \cite[Theorems 5 and 6]{our} implies that
$\exp(H_2\cap\Omega^{\varepsilon}_{2n}(q))=\exp(H_2)$, if $a+b>1$, and
$\exp(H_2\cap\Omega^{\varepsilon}_{2n}(q))=\exp(H_2)/2$, if $a+b=1$. Assume that $Z\neq 1$. Since group $H_2\cap
\Omega^\varepsilon_{2n}(q)$ intersects with $Z$ trivially, we have $\eta(H)=\eta(\tilde H)$ for $a+b>1$. Let $a+b=1$.
Then $H_2$ is cyclic. Denote by $t$ a generator of $H_2$. We have $Z(H\cap \Omega^\varepsilon_{2n}(q))=\langle
t_0^{|t_0|/2}\rangle\times\langle t^2\rangle$, and $Z=\langle t_0^{|t_0|/2}t^{|t|/2}\rangle$. Since $Z\subseteq H\cap
\Omega^\varepsilon_{2n}(q)$, the number $|t|/2$ is even. Therefore, the center of $H\cap \Omega^\varepsilon_{2n}(q))$
can be presented as the direct product $\langle t_0^{|t_0|/2}t^{|t|/2}\rangle\times\langle t^2\rangle$. This implies
that $Z(\tilde H)\simeq \langle t^2\rangle$ and, therefore, $\eta(\tilde H)=\eta(H)/2$.

Let $|t_0|/2$ be odd. Then $Z(H\cap \Omega^{\varepsilon}_{2n}(q))$ is conjugate in $\overline H$ to
$$\prod_{t_i\neq t_0}\langle t_0^{|t_0|/2}t_i\rangle\simeq H_2.$$
Therefore, $\eta(H\cap\Omega^\varepsilon_{2n}(q))=\eta(H)$. Assume that $Z\neq 1$. Lemma \ref{l} implies that if
$a+b>1$, then $\eta(\tilde H)=\eta(H)$. If $a+b=1$, then $\eta(\tilde H)=\eta(H)/2$. The first statement is proved.

$2)$  Let $T_1$ be a maximal torus of $H_1$ isomorphic to a cyclic group of order $(q^2-1)$. Then the subgroup $T$ of
$H$, equal to $T_1\times H_2$, is a maximal torus of $SO^\varepsilon_{2n}(q)$. Define a pair of partitions
$\alpha=(2,\alpha_1\dots,\alpha_2)$ and $\beta=(\beta_1,\dots,\beta_b)$. Then structure of $T$ is determined by the
pair of partitions $\alpha$ and $\beta$, and by Lemma \ref{l:toriOmega2n} the torus $T\cap\Omega^\varepsilon_{2n}(q)$
is conjugate in $\overline H$ to $\{t_1^{k_1}t_2^{k_2}\dots t_{a+b+1}^{k_{a+b+1}}| k_1+k_2+\dots+k_{a+b+1}\text{ is
even}\}$, where elements $t_1,\dots,t_{a+b+1}$ are defined before the statement of Lemma \ref{l:toriOmega2n}. We have
$Z(H)=\langle t_1^{q+\varepsilon_11}\rangle\times H_2$. Since the number $q+\varepsilon_11$ is even, we obtain that
$Z(H\cap\Omega^\varepsilon_{2n}(q))=\langle t_1^{q+\varepsilon_11}\rangle\times (H_2\cap\Omega^\varepsilon_{2n}(q))$.
Now, \cite[Theorems 5 and 6]{our} implies that $\exp(H_2\cap\Omega^{\varepsilon}_{2n}(q))=\exp(H_2)$, if $a+b>1$, and
$\exp(H_2\cap\Omega^{\varepsilon}_{2n}(q))=\exp(H_2)/2$, if $a+b=1$.

Assume that $Z\neq 1$. By \cite[Theorems 5 and 6]{our} there exist elements $h_1$, $h_2,\dots, h_{a+b}$ in  $H_2$ such
that $h_1\not\in\Omega^\varepsilon_{2n}(q)$,
$|h_1|_{\{2\}}\leqslant|h_2|_{\{2\}}\leqslant\dots\leqslant|h_{a+b}|_{\{2\}}$ and $H_2=\langle h_1\rangle\times\langle
h_2\rangle\times\dots\times\langle h_{a+b}\rangle$. Moreover, if the number
$(|h_1|_{\{2\}}+|h_2|_{\{2\}}+\dots+|h_{a+b}|_{\{2\}})/|h_1|_{\{2\}}$ is even, then the projection $Z$ on $H_2$ is
contained in the subgroup $\langle h_2\rangle$, and if it is odd, then the projection is contained in $\langle
h_1^2\rangle$. We have $Z\leqslant\langle t_0^{q+\varepsilon_11}\rangle\times\langle h_1^2\rangle$ or
$Z\leqslant\langle t_0^{q+\varepsilon_11}\rangle\times\langle h_2\rangle$. In both cases by applying Lemma \ref{l} we
obtain that $\eta(\tilde H)=\eta(H\cap\Omega^\varepsilon_{2n}(q))$. Lemma \ref{l:mu(SO2n)} is proved.
\end{prf}

\begin{theorem}\label{t:mum(omega2n)} Let $G=\Omega^\varepsilon_{2n}(q)$, where $n\geqslant 4$, $\varepsilon\in\{+,-\}$
and $q$ is a power of an odd prime $p$. For  $k\geqslant 1$ put $n(k)=(p^{k-1}+3)/2$. Let $\nu(G)$ to consist of the
following numbers:
\begin{itemize}

\item[$1)$] $p^k[d_k,\frac{q^{n-n(k)}\pm1}{d_k}]$, где $d_k=\frac{(4, q^{n(k)}-\varepsilon1)}{2}$,
for every $k$ such that $n(k)<n$;

\item[$2)$] $p^k[q^{n_1}+1, q^{n_2}+1, \dots, q^{n_l}+1, q^{n_{l+1}}-1, q^{n_{l+2}}-1,\dots, q^{n_s}-1]$
for all $s>1$ and $n_1,n_2,\dots,n_s>0$ such that $n(k)+n_1+n_2+\dots+n_s=n$;

\item[$3)$] $p[q\pm1, q^{n_1}+1, q^{n_2}+1, \dots, q^{n_l}+1, q^{n_{l+1}}-1, q^{n_{l+2}}-1,\dots, q^{n_s}-1]$
for all $s>1$, $l$ is even, if $\varepsilon=+$, and odd, if $\varepsilon=-$, and $n_1,n_2,\dots,n_s>0$ such that
$2+n_1+n_2+\dots+n_s=n$;

\item[$4)$] $p[q\pm1, \frac{q^{n-2}-\varepsilon1}{2}]$;

\item[$5)$] $2p^k$, if $n=n(k)$ and $(4, q^n-\varepsilon1)=4$;

\item[$6)$] $p(q^2\pm1)$, if $n=4$ and $\varepsilon=+$.

\item[$7)$] $9(q\pm1)$, if $n=4$, $p=3$ and $\varepsilon=+$.
\end{itemize}
Then $\mu_m(G)\subseteq \nu(G) \subseteq \omega(G)$.
\end{theorem}

\begin{theorem}\label{t:mum(Pomega2n)} Let $G=P\Omega^\varepsilon_{2n}(q)$, where $n\geqslant 4$, $\varepsilon\in\{+,-\}$
and $q$ is a power of an odd prime $p$, and $(4,q^n-\varepsilon1)=4$. For $k\geqslant 1$ put $n(k)=(p^{k-1}+3)/2$. Let
$\nu(G)$ to consist of the following numbers:
\begin{itemize}

\item[$1)$] $p^k\frac{q^{n-n(k)}\pm1}{2}$
for every $k$ such that $n(k)<n$;

\item[$2)$] $p^k[q^{n_1}+1, q^{n_2}+1, \dots, q^{n_l}+1, q^{n_{l+1}}-1, q^{n_{l+2}}-1,\dots, q^{n_s}-1]$
for all $s>1$ and $n_1,n_2,\dots,n_s>0$ such that $n(k)+n_1+n_2+\dots+n_s=n$;

\item[$3)$] $p[q\pm1, q^{n_1}+1, q^{n_2}+1, \dots, q^{n_l}+1, q^{n_{l+1}}-1, q^{n_{l+2}}-1,\dots, q^{n_s}-1]$
for all $s>1$, $l$ is even, if $\varepsilon=+$, and odd, if $\varepsilon=-$, and $n_1,n_2,\dots,n_s>0$ such that
$2+n_1+n_2+\dots+n_s=n$;

\item[$4)$] $p[q\pm1, \frac{q^{n-2}-\varepsilon1}{2}]$.
\end{itemize}
Then $\mu_m(G)\subseteq \nu(G) \subseteq \omega(G)$.
\end{theorem}

We will prove Theorems \ref{t:mum(omega2n)} and \ref{t:mum(Pomega2n)} simultaneously.

\begin{prf} Let $\varepsilon=+$.

Let us prove that  $\nu(G)\subseteq \omega(G)$. Notice that the number $2p^k$ in item $5)$ of Theorem
\ref{t:mum(omega2n)} arises as $\eta(\Omega^\varepsilon_{2n}(q))$, since the condition $(4, q^n-\varepsilon1)=4$ is
equivalent to the condition $Z(\Omega^\varepsilon_{2n}(q))\neq 1$. If $n\neq n(k)$ for every $k$, then Lemma
\ref{l:mu(SO2n)} implies that $\eta(G)$ divides $\eta(H)$ for some proper reductive subgroup $H$. Further all
considered reductive subgroup are proper. Let $n=4$. Then $SO^+_8(q)$ contains subgroups $H_1$ and $H_2$ isomorphic to
$GL_2(q^2)$ and $GU_2(q^2)$ respectively. The centers of $H_1$ and $H_2$ are contained in maximal tori of these groups
isomorphic to a cyclic group of order $(q^{4}-1)$. Since numbers $(q^4-1)/(q^2-1)$ and $(q^4-1)/(q^2+1)$ are even, in
the same way as in the proof  of Lemma \ref{l:mu(SO2n)} we obtain that $Z(H_1)$ and $Z(H_2)$ are contained in
$\Omega^+_{8}(q)$. Thus, we deduce that $\eta(H_1)$ and $\eta(H_2)$ lie in $\omega(\Omega^+_{8}(q))$. Let $n=4$ and
$p=3$. Then $SO^+_8(q)$ contains subgroups $H_1$ and $H_2$ isomorphic to $GL_4(q)$ and $GU_4(q)$ respectively. In this
case generators of the centers of these subgroups are again even powers of generators of cyclic tori of orders $q^4-1$.
Thus, we deduce that $\eta(H_1)$ and $\eta(H_2)$ lie in $\omega(\Omega^+_{8}(q))$. The remaining inclusions follow from
Lemma \ref{l:mu(SO2n)}.

Let us prove the inclusion $\mu_m(G)\subseteq \nu(G)$. Let $H$ be a proper reductive subgroup of maximal rank of
$SO^+_{2n}(q)$, which is not a maximal torus. Denote by $H_1$ the group $H\cap\Omega^+_{2n}(q)$, and by $H_2$ the image
of $H_1$ in $P\Omega^+_{2n}(q)$. Let structure of $H$ be determined by partitions $\xi^{(i)}$ and $\zeta^{(i)}$, where
$1\leqslant i\leqslant a$, $\theta^{(i)}$ and $\upsilon^{(i)}$, where $2\leqslant i\leqslant b$. Since $H$ is not a
maximal torus, either $a>1$ or $b>1$. By formula (\ref{zorth}) the center of $H$ is isomorphic to the direct product
$$\prod_{i,j}(q^{\xi^{(i)}_j}-1)\times\prod_{i,j}(q^{\zeta^{(i)}_j}+1)\times
\prod_{i,j}Z(H_{ij}^\theta)\times\prod_{i,j}Z(H_{ij}^\upsilon).$$ Put $e=\exp(
\prod_{i,j}Z(H_{ij}^\theta)\times\prod_{i,j}Z(H_{ij}^\upsilon))$. We have
$$\eta(H)=p^k[\underset{i,j}{\operatorname{lcm}}\{q^{\xi^{(i)}_j}-1\},\underset{i,j}{\operatorname{lcm}}
\{q^{\zeta^{(i)}_j}+1\}, e],$$ where $k$ is such that $p^k$ is the least power of $p$ which is greater than
$\operatorname{max}\{a-1, 2b-3\}$, i.\,e. the greatest power of $p$ lying in $\omega(H)$. We have
$\eta(H_1)=\eta(H_1)/c_1$ and $\eta(H_2)=\eta(H)/c_2$, where $c_1\in\{1,2\}$ and $c_2\in\{1,2,4\}$ and specific values
of $c_1$ and $c_2$ depend on structure of $H$. Put $n_0=(p^{k-1}+3)/2$.

Put
$$x=|\xi^{(1)}|+|\zeta^{(1)}|+|\xi^{(2)}|+|\zeta^{(2)}|+\dots+|\xi^{(a)}|+|\zeta^{(a)}|+n_0.$$
By word for word repetition of arguments from the proof of Theorem \ref{t:so2n} we obtain that $x>n$ implies that group
$H$ is isomorphic to
$$GL_2(q)\times\prod_{j}(q^{\xi^{(1)}_j}-1)\times\prod_{j}(q^{\zeta^{(1)}_j}+1),$$ or
$$GU_2(q)\times\prod_{j}(q^{\xi^{(1)}_j}-1)\times\prod_{j}(q^{\zeta^{(1)}_j}+1).$$ In both cases numbers
$\eta(H_1)$ and $\eta(H_2)$ lie in corresponding sets $\nu(G)$.

Let $x\leqslant n$. If the sum of lengths of partitions $\zeta^{(i)}$ over all $i$ is even, then $SO^+_{2n}(q)$
contains a subgroup $K$ isomorphic to the direct product

$$SO^+_{2n_0}(q)\times\prod_{i,j}(q^{\xi^{(i)}_j}-1)\times\prod_{i,j}(q^{\zeta^{(i)}_j}+1)\times
\underbrace{(q-1)\times(q-1)\times\dots\times(q-1)}_{n-x},$$ if the sum is odd, then $SO^+_{2n}(q)$ contains a subgroup
$K$ isomorphic to the direct product

$$SO^-_{2n_0}(q)\times\prod_{i,j}(q^{\xi^{(i)}_j}-1)\times\prod_{i,j}(q^{\zeta^{(i)}_j}+1)
\times \underbrace{(q-1)\times(q-1)\times\dots\times(q-1)}_{n-x}.$$ We have
$$\eta(K)=p^k[\underset{i,j}{\operatorname{lcm}}\{q^{\xi^{(i)}_j}-1\},\underset{i,j}{\operatorname{lcm}}
\{q^{\zeta^{(i)}_j}+1\}, b],$$ where $b=q-1$, if $n-x>0$, and $b=1$, if $n-x=0$. Denote by $K_1$ the intersection
$K\cap\Omega^+_{2n}(q)$ and by $K_2$ the image of $K_1$ in $P\Omega^+_{2n}(q)$. If the sum of lengths of partitions
$\xi^{(i)}$, $\zeta^{(i)}$ over all $i$ and of number $(n-x)$ is greater than $1$, then Lemma \ref{l:mu(SO2n)} implies
that $\eta(K_1)=\eta(K_2)=\eta(K)$. Thus, $\eta(K)$ is divided by $\eta(H)$ and, therefore, by $\eta(H_1)$ and
$\eta(H_2)$, and lies in $\nu(G)$. Assume that the sum of lengths of partitions $\xi^{(i)}$, $\zeta^{(i)}$ over all $i$
and of number $(n-x)$ equals $1$. Then one of the following conditions holds:

1) $n-x=0$ and the sum of lengths of partitions $\xi^{(i)}$ and $\zeta^{(i)}$ over all $i$ is equal to $1$;

2) $n-x=1$ and $a=0$.

Consider the first case. Recall that
$$n=|\xi^{(1)}|+|\zeta^{(1)}|+2(|\xi^{(2)}|+|\zeta^{(2)}|)+\dots+a(|\xi^{(a)}|+|\zeta^{(a)}|)+$$
$$+2(|\theta^{(2)}|+|\upsilon^{(2)}|)+\dots+b(|\theta^{(b)}|+
|\upsilon^{(b)}|).$$

By substituting in equality $n=x$ the expressions for $n$ and $x$, we obtain the following equality
$$n_0=|\xi^{(2)}|+|\zeta^{(2)}|+\dots+(a-1)(|\xi^{(a)}|+|\zeta^{(a)}|)
+2(|\theta^{(2)}|+|\upsilon^{(2)}|)+\dots+b(|\theta^{(b)}|+ |\upsilon^{(b)}|).$$ Therefore,
$n_0\geqslant\max\{a-1,b\}$. We have $n_0\leqslant\max\{a/2+1,b\}$.

Assume that $\max\{(a+3)/2, b\}=b$. Then $n_0=b$, $|\theta^{(b)}|+ |\upsilon^{(b)}|=1$, $|\theta^{(i)}|+
|\upsilon^{(i)}|=0$ for $i<b$ and $a=1$. This implies that $H$ is isomorphic to
$$SO^{\varepsilon_1}_{2n_0}(q)\times\prod_{i,j}(q^{\xi^{(i)}_j}-1)\times\prod_{i,j}(q^{\zeta^{(i)}_j}+1),$$
where $\varepsilon_1=+$, if $|\theta^{(b)}|=1$, and $\varepsilon_1=-$, if $|\upsilon^{(b)}|=1$. Lemma \ref{l:mu(SO2n)}
implies that in this case numbers $\eta(H_1)$ and $\eta(H_2)$ lie in corresponding sets $\nu(G)$.

Assume that $\max\{a/2+1, b\}=a/2+1$. Then $2\leqslant a\leqslant 4$. Let $a=4$. Then $n_0=3$, $b=0$, $n=4$ и $p=3$.
Thus, either $H\simeq GL_4(q)$ or $H\simeq GU_4(q)$. As it was shown above, in this case $\eta(H_1)$ equals either
$9(q-1)$ or $9(q+1)$ and, therefore, lies in $\nu(\Omega^+_8(q))$. The center of $H_2$ in this case is isomorphic to a
cyclic group of order $(q\pm1)/2$. Thus, $\eta(H_2)=3^2(q\pm1)/2$. We obtain that $\eta(H_2)$ lies in
$\nu(P\Omega^+_8(q))$. Notice that this case is impossible for $\varepsilon=-$, since $\Omega^-_8(q)$ does not contain
subgroups isomorphic to $GL_4(q)$ or $GU_4(q)$. Let $a=3$. Then $n_0=2$, $b=0$ и $n=3$. The last contradicts with the
hypothesis of the theorem. Let $a=2$. Then $n_0=2$. We have
$n_0=2=(2-1)(|\xi^{(2)}|+|\zeta^{(2)}|)+2(|\theta^{(2)}|+|\upsilon^{(2)}|)+\dots+b(|\theta^{(b)}|+ |\upsilon^{(b)}|)$.
This implies that $b=0$ and $|\xi^{(2)}|+|\zeta^{(2)}|=2$. Since the sum of lengths of partitions $\xi^{(i)}$ and
$\zeta^{(i)}$ over all $i$ equals $1$, Proposition \ref{p:so+-fin} implies that $H$ is isomorphic to $GL_2(q^2)$, if
$|\xi^{(2)}|=2$, and isomorphic to $GU_2(q^2)$, if $|\zeta^{(2)}|=2$. We have $n=4$. Then, as it was shown above,
$\eta(H_1)$ is equal to $p(q^2-1)$, if $H\simeq GL_2(q^2)$, equal to $p(q^2+1)$, if $H\simeq GГ_2(q^2)$, and lies in
$\nu(\Omega^+_8(q))$. The center of $H_2$ in this case is isomorphic to a cyclic group of order $(q^2\pm1)/2$. Thus,
$\eta(H_2)=p(q^2\pm1)/2$. We deduce that $\eta(H_2)$ lies in $\nu(P\Omega^+_8(q))$. Notice that this case is also
impossible for $\varepsilon=-$, since $\Omega^-_8(q)$ does not contain subgroups isomorphic to $GL_2(q^2)$ or
$GU_2(q^2)$.

Assume that $n-x=1$ and $a=0$. We have $$n_0=2(|\theta^{(2)}|+|\upsilon^{(2)}|)+\dots+b(|\theta^{(b)}|+
|\upsilon^{(b)}|)-1.$$ Since $n_0\leqslant b$, we deduce that $n_0=b$ and $H=\Omega^+_{2n}(q)$.

The proof for the case $\varepsilon=-$ can be done in the same way.

The theorems are proved.
\end{prf}

As a corollary of Lemma \ref{l:pperoid}, \cite[Theorems 5 and 6]{our} and Theorems \ref{t:mum(omega2n)} and
\ref{t:mum(Pomega2n)} we obtain the following statements.

\begin{cor} Let $G=\Omega_{2n}^\varepsilon(q)$, where $n\geqslant 4$, $\varepsilon\in\{+,-\}$ and $q$ is a power of an odd
prime $p$. For $k\geqslant 1$ put $n(k)=(p^{k-1}+3)/2$. Then $\omega(G)$ consists of all divisors of the following
numbers:
\begin{itemize}
\item[$1)$] $\frac{q^{n}-\varepsilon1}{2}$;

\item[$2)$] $[q^{n_1}+1, q^{n_2}+1, \dots, q^{n_l}+1, q^{n_{l+1}}-1, q^{n_{l+2}}-1,\dots, q^{n_s}-1]$
for all $s>1$, $l$ is even, if $\varepsilon=+$, and odd, if $\varepsilon=-$, and $n_1,n_2,\dots,n_s>0$ such that
$n_1+n_2+\dots+n_s=n$;

\item[$3)$] $p^k[d_k,\frac{q^{n-n(k)}\pm1}{d_k}]$, where $d_k=\frac{(4, q^{n(k)}-\varepsilon1)}{2}$,
for every $k$ such that $n(k)<n$;

\item[$4)$] $p^k[q^{n_1}+1, q^{n_2}+1, \dots, q^{n_l}+1, q^{n_{l+1}}-1, q^{n_{l+2}}-1,\dots, q^{n_s}-1]$
for all $s>1$ and $n_1,n_2,\dots,n_s>0$ such that $n(k)+n_1+n_2+\dots+n_s=n$;

\item[$5)$] $p[q\pm1, q^{n_1}+1, q^{n_2}+1, \dots, q^{n_l}+1, q^{n_{l+1}}-1, q^{n_{l+2}}-1,\dots, q^{n_s}-1]$
for all $s>1$, $l$ is even, if $\varepsilon=+$, and odd, if $\varepsilon=-$, and $n_1,n_2,\dots,n_s>0$ such that
$2+n_1+n_2+\dots+n_s=n$;

\item[$6)$] $p[q\pm1, \frac{q^{n-2}-\varepsilon1}{2}]$;

\item[$7)$] $dp^k$, where $d=(4, q^n-\varepsilon1)/2$, if $n=n(k)$ for some $k$;

\item[$8)$] $p(q^2\pm1)$, if $n=4$ and $\varepsilon=+$.

\item[$9)$] $9(q\pm1)$, if $n=4$, $p=3$ and $\varepsilon=+$.
\end{itemize}
\end{cor}

\begin{cor}  Let $G=P\Omega_{2n}^\varepsilon(q)$, where $n\geqslant 4$, $\varepsilon\in\{+,-\}$, $q$ is a
power of an odd prime $p$, and $(4,q^n-\varepsilon1)=4$. For $k\geqslant 1$ put $n(k)=(p^{k-1}+3)/2$. Then $\omega(G)$
consists of all divisor of the following numbers:
\begin{itemize}
\item[$1)$] $\frac{q^{n}-\varepsilon1}{4}$;

\item[$2)$] $\frac{[q^{n_1}-\varepsilon_11,q^{n_2}-\varepsilon\varepsilon_11]}{d}$, where $n_1+n_2=n$,
$\varepsilon_1\in\{+,-\}$, $d=2$, if
$(q^{n_1}-\varepsilon_11)_{\{2\}}=(q^{n_2}-\varepsilon\varepsilon_11)_{\{2\}}$, and $d=1$
otherwise;

\item[$3)$] $[q^{n_1}+1, q^{n_2}+1, \dots, q^{n_l}+1, q^{n_{l+1}}-1, q^{n_{l+2}}-1,\dots, q^{n_s}-1]$
for all $s>2$, $l$ is even, if $\varepsilon=+$, and odd, if $\varepsilon=-$, and $n_1,n_2,\dots,n_s>0$ such that
$n_1+n_2+\dots+n_s=n$;

\item[$4)$] $p^k\frac{q^{n-n(k)}\pm1}{2}$
for every $k$ such that $n(k)<n$;

\item[$5)$] $p^k[q^{n_1}+1, q^{n_2}+1, \dots, q^{n_l}+1, q^{n_{l+1}}-1, q^{n_{l+2}}-1,\dots, q^{n_s}-1]$
for all $s>1$ and $n_1,n_2,\dots,n_s>0$ such that $n(k)+n_1+n_2+\dots+n_s=n$;

\item[$6)$] $p[q\pm1, q^{n_1}+1, q^{n_2}+1, \dots, q^{n_l}+1, q^{n_{l+1}}-1, q^{n_{l+2}}-1,\dots, q^{n_s}-1]$
for all $s>1$, $l$ is even, if $\varepsilon=+$, and odd, if $\varepsilon=-$, and $n_1,n_2,\dots,n_s>0$ such that
$2+n_1+n_2+\dots+n_s=n$;

\item[$7)$] $p[q\pm1, \frac{q^{n-2}-\varepsilon1}{2}]$;

\item[$8)$] $p^k$, if $n=n(k)$ for some $k$.
\end{itemize}
\end{cor}

\end{document}